\documentclass{article}
\usepackage[utf8]{inputenc}
\usepackage[margin=3.5cm]{geometry}

\usepackage{amsmath}
\usepackage{amsfonts}
\usepackage{amssymb}
\usepackage{graphicx}
\usepackage{epstopdf}
\usepackage{subfigure}

\usepackage{color}

\newcommand{\e}{{\rm e}}        
\def\R{\mathbb{ R}}             

\renewcommand{\d}{{\rm d}}      
\def\dW{{\rm d}W}               
\def\dt{{\rm d}t}

\def\ds{{\rm d}s}

\def\1{{\mathds{1}}}            

\title{A stochastic $\theta$-SEIHRD model: adding randomness to the COVID-19 spread}
\author{Álvaro Leitao and Carlos Vázquez \\
Department of Mathematics and CITIC, Universidade da Coruña, Spain}
\date{\today}

\begin{document}

\maketitle

\begin{abstract}
In this article we mainly extend the deterministic model developed in \cite{ivorra2020} to a stochastic setting. More precisely, we incorporated randomness in some coefficients by assuming that they follow a prescribed stochastic dynamics. In this way, the model variables are now represented by stochastic process, that can be simulated by appropriately solve the system of stochastic differential equations. Thus, the model becomes more complete and flexible than the deterministic analogous, as it incorporates additional uncertainties which are present in more realistic situations. In particular, confidence intervals for the main variables and worst case scenarios can be computed.
\end{abstract}

\section{Introduction}
    
    The coronavirus disease 2019, renamed as COVID-19, is an infectious disease produced by the virus SARS-CoV-2, which was declared pandemic by the World Health Organization (WHO) in March 11, 2020. This disease has posed many novel scientific challenges, ranging from contagious patterns, medical treatments or vaccine developments, to data analytics, spread modelling or evolution forecast. The research on all of these topics has been intensive in the last few months, as so will be in the near future. Thus, many disciplines from interconnected fields, like bio-sciences, mathematical and statistical modelling or artificial intelligence, will be involved and they will play an important role to rapidly overcome this exceptional emergency situation.
    
    From the mathematical modelling point of view, epidemiological compartmental models have been often used to understand and analyze the behaviour of the contagious diseases, with the article \cite{kermack1927} being the pioneering work. These models provide useful tools to make predictions about the future evolution of the epidemic and control its propagation. In the literature of epidemic modelling, many examples of general-purpose models have been proposed (see \cite{hethcote2000}, for a review), each of them accounting for some specific characteristics of the diseases, like the well-known examples SIR (Susceptible-Infected-Recovered), SEIR (Susceptible-Exposed-Infected-Recovered) or SIS (Susceptible-Infected-Susceptible) models. 
    
    Typically, the compartmental models are originally formulated following a deterministic approach, i.e., in terms of a system of \emph{Ordinary Differential Equations} (ODEs), although usually they are readily extended to an stochastic version. There are two common approaches to include stochasticity into a deterministic model, relying on either the \emph{Continuous Time Markov Chain} (CTMC) or \emph{Stochastic Differential Equations} (SDEs), see \cite{allen2017} and the references therein. Here, we will focus on the latter. The stochastic models allow to capture many kinds of circumstances including uncertainty that may influence the compartments dynamics: behavioural effects, public interventions, seasonal patterns, environmental factors, etc. Furthermore, the solution of a deterministic model is given by a set of functions of time uniquely dependent on the initial data, while the solution of the stochastic model is a set of stochastic processes, containing much more information than the deterministic analogous. In fact, at each time instant we can exploit the information provided by the probability distribution associated to the underlying random variable. An statistical analysis can therefore be performed, producing useful outcomes like quartiles or worst case scenarios.
    
    In this work we present an stochastic extension of the deterministic compartmental model in \cite{ivorra2020}, that has been proposed as an \emph{ad hoc} model for the COVID-19 disease. Unlike to the classical models (SIR, SEIR, SIS), this model is adapted to the specific characteristics of the COVID-19, taking into account, besides the usual factors, the undetected infectious cases, the hospitalized population or the deaths, for example. Thus, the so-called $\theta$-SEIHRD model proposed in \cite{ivorra2020} was developed as a very general and sophisticated model (based on a previously introduced model \cite{ivorra2015}) in order to be able to study the spread of the disease worldwide. In the same work, the authors proposed a simplified version which reduces the mathematical complexity towards a more tractable model, but yet preserving the ability to capture the most relevant features. Other recent studies addressing the mathematical or statistical modelling in the context of the COVID-19 include \cite{kucharski2020, roosa2020, platen2020, lipton2020}, for example. As the literature focusing on COVID-19 modelling is increasing a lot during last months and weeks, it results quite difficult for the authors to select a list of the more relevant papers more recently arising in the topic, surely some of them being developed in parallel to this contribution.
    
    As mentioned before, we will follow the SDE approach where the randomness is typically achieved by incorporating a \emph{Brownian motion}, also known as \emph{Wiener process}, to the ODEs of the deterministic system. There are two common ways of addressing this stochastic extension. On the one hand, an arbitrary random noise can be added to some of the equations in the ODE system, thus transforming them into SDEs. On the other hand, one (or more) of the existing model parameters can be perturbed, meaning that it (they) becomes a random variable or a stochastic process. Although both approaches have been well investigated, recent examples are \cite{han2018,li2018} for the random noise approach and \cite{gray2011, dureau2013, cai2019} for the parameter perturbation, here we prefer the second alternative. When adding random noise to an equation, it often comes with an extra parameter to control the level of volatility, which might have no biological interpretation. However, in practice, the uncertainty will have impact on a particular model component, typically represented by a model parameter. Therefore, a randomly perturbed parameter can be reasonably explained in terms of the variability produced by the source of the considered uncertainty.
    
    In this work, we will consider a randomly perturbed disease contact rates. The approach proposed here is however conceptually different than the ones typically found in the literature, aiming to deal with some of the observed inconsistencies in the COVID-19 context. As a natural choice, many stochastic SIR-like models rely on a normally distributed perturbation for the disease contact rate parameters (see \cite{gray2011, cai2019}, for example). However, this approach allows negative values for the rates, which is biologically non-sensical, potentially appearing when the rate is close to zero, as it is seen in practice for the COVID-19 \cite{li2020}. This issue is overcome in \cite{dureau2013}, by employing an exponential Ornstein-Uhlenbeck (OU) process to model the contact rate. Although this process indeed ensures positiveness, it presents another undesirable feature: an increasing variance in time. When analysing the disease transmission dynamics, there is not an objective evidence of more variability in the disease contact rates in the long-term \cite{whoreports2020}. Actually, one can expect more control of the disease spread patterns over the course of the epidemic. For all the reasons mentioned above, we propose to model the contact disease rates by means of the so-called \emph{CIR process}, named after its authors Cox, Ingersoll and Ross \cite{cox1985}. The CIR process is widely employed to simulate the evolution of interest rates in the quantitative finance framework. In some sense, the interest rates in computational finance and the disease contact rates in epidemiology present a rather similar behaviour: positiveness, controlled variability and long-term stability.
    
    The article is organized as follows. In Section 2 we introduce the proposed mathematical model. Section 3 is describes the numerical methods for solving the model. Section 4 contains the obtained results and their numerical and statistical analysis. Finally, Section 5 points out some conclusions and possible future research lines.
    

\section{A stochastic $\theta$-SEIHRD model of the COVID-19}

    As mentioned in the previous introduction, we base our approach in the model in \cite{ivorra2020}, which provides an advanced extension of the classical SIR-like models that adds to the usual compartments in the SEIR model the very specific ones for COVID-19: undetected infectious, hospitalized and deaths. In \cite{ivorra2020}, the authors first presented a very sophisticated, well motivated and general model, with a significant amount of free parameters allowing an enormous flexibility, including several territories. In order to make the model more practically tractable, the authors also proposed a simplified version, by imposing several assumptions: single territory, regular natality/mortality is neglected, and no imported/exported cases. Furthermore, some predefined forms for the open parameters of the model are adopted.
    
    Thus, the simplified $\theta$-SEIHRD model proposed in \cite{ivorra2020} is therefore given by
    \begin{equation}\label{eq:thetaSEIHRD_model}
    \begin{aligned}
        \frac{\d S}{\d t}(t) &= -\frac{S(t)}{N}\left(m_E(t) \beta_E E(t) + m_I(t) \beta_I I(t) + m_{I_u}(t) \beta_{I_u}(\theta(t)) I_u(t)\right) \\ 
        & \quad -\frac{S(t)}{N}\left(m_{H_R}(t) \beta_{H_R}(t) H_R(t) + m_{H_D}(t) \beta_{H_D}(t) H_D(t)\right), \\
        \frac{\d E}{\d t}(t) &= \frac{S(t)}{N}\left(m_E(t) \beta_E E(t) + m_I(t) \beta_I I(t) + m_{I_u}(t) \beta_{I_u}(\theta(t)) I_u(t)\right) \\
        &+ \frac{S(t)}{N}\left(m_{H_R}(t) \beta_{H_R}(t) H_R(t) + m_{H_D}(t) \beta_{H_D}(t) H_D(t)\right) - \gamma_E E(t), \\
        \frac{\d I}{\d t}(t) &= \gamma_E E(t) - \gamma_I(t) I(t), \\
        \frac{\d I_u}{\d t}(t) &= (1 - \theta(t))\gamma_I(t) I(t) - \gamma_{I_u}(t) I_u(t), \\
        \frac{\d H_R}{\d t}(t) &= \theta(t)\left(1 - \frac{\omega(t)}{\theta(t)}\right)\gamma_I(t) I(t) - \gamma_{H_R}(t) H_R(t), \\
        \frac{\d H_D}{\d t}(t) &= \omega(t)\gamma_I(t) I(t) - \gamma_{H_D}(t) H_D(t), \\
        \frac{\d R_d}{\d t}(t) &= \gamma_{H_R}(t) H_R(t), \\
        \frac{\d R_u}{\d t}(t) &= \gamma_{I_u}(t) I_u(t), \\
        \frac{\d D}{\d t}(t) &= \gamma_{H_D}(t) H_D(t). 
    \end{aligned}
    \end{equation}
    Note that the last three equations of the system (\ref{eq:thetaSEIHRD_model}) are uncoupled with previous ones, so that the expression of their solution can be obtained by
    \begin{equation*}
    \begin{aligned}
        R_d(t) &= R_d(t_0) + \int_{t_0}^t \gamma_{H_R}(s)H_R(s)\ds, \\
        R_u(t) &= R_u(t_0) + \int_{t_0}^t \gamma_{I_u}(s)I_u(s)\ds, \\
        D(t) &= D(t_0) + \int_{t_0}^t \gamma_{H_D}(s)H_D(s)\ds,
    \end{aligned}
    \end{equation*}
    once the first six ones have been solved in a first stage. Thus, only the first six coupled equations need to be solved. The system (\ref{eq:thetaSEIHRD_model}) is completed with the initial data $S(t_0)$, $E(t_0)$, $I(t_0)$, $I_u(t_0)$, $H_R(t_0)$, $H_D(t_0)$, $R_d(t_0)$, $R_u(t_0)$ and $D(t_0)$.
    
    The model includes a number of open parameters, the model parameters, which are often determined either from the available data/literature or by calibration. In the following, a detailed explanation of each parameter is provided:
    \begin{itemize}
        \item Efficiency of the control measures $m_E, m_I, m_{I_u}, m_{H_R}, m_{H_D} \in [0, 1] (\%)$. Here, only one control measure is assumed (mobility restrictions, for example), implemented at date $\lambda_1$ and represented by a decaying function
            \begin{equation*}
                m_E(t) = m_I(t) = m_{I_u}(t) = m_{H_R}(t) = m_{H_D}(t) = \left\{\begin{aligned}
                                                                    &1, & \text{ if } t \in [0, \lambda_1], \\
                                                                    &\exp\left(-\kappa_1(t - \lambda_1)\right), & \text{ if } t \in [\lambda_1, T],
                                                                \end{aligned}
                                                        \right.
            \end{equation*}
            with the parameter $\kappa_1 \in [0, 0.2]$ entering in the calibration procedure. The generalization to more/individualized control measures is straightforward.
        
        \item The fatality rate $\omega(t) \in [\underline{\omega}, \overline{\omega}] \subset [0,1]$, for which the following form is proposed
            \begin{equation*}
                \omega(t) = m_I(t)\overline{\omega} + (1 - m_I(t))\underline{\omega},
            \end{equation*}
            with $\underline{\omega}$ and $\overline{\omega}$ being the fatality rate limits with and without control measures, respectively. As $\overline{\omega} \geq \underline{\omega}$, its value is defined in terms of a displacement w.r.t the lower limit, i.e. $\overline{\omega} = \underline{\omega} + \delta_\omega$, whose values are calibrated.
        
        \item The fraction of infected individuals which are detected and documented by the authorities, $\theta$. Assuming that all deaths are detected i.e. $\theta \in [\overline{\omega}, 1]$ and it changes in time,
            \begin{equation*}
                \theta(t) = \left\{\begin{aligned}
                                    &\underline{\theta}, & \text{ if } t \in [t, \lambda_1], \\
                                    &\text{linear continuous}, & \text{ if } t \in [\lambda_1, \lambda_2], \\
                                    &\overline{\theta}, & \text{ if } t \in [\lambda_2, T],
                                    \end{aligned}
                            \right.
            \end{equation*}
            with $\underline{\theta}, \overline{\theta}, \lambda_1, \lambda_2$ inferred from the data.
        
        \item The disease contact rates $\beta_E, \beta_I, \beta_{I_u}, \beta_{H_R}, \beta_{H_D} \in \R^+$. The parameter $\beta_I$ is assumed to be given, calibrated to the available data, and a relation between it and the rest of the contact rates exists,
            \begin{equation}\label{eq:betas}
                \beta_E = C_E\beta_I, \quad
                \beta_{I_u}(t) = \underline{\beta}_I + \frac{\beta_I - \underline{\beta}_I}{1 - \omega(t)}(1 - \theta(t)), \quad
                \beta_{H_R} = \beta_{H_D} = C_H(t)\beta_I,
            \end{equation}
            where $\underline{\beta}_I = C_u\beta_I$, with $C_E$, $C_H(t)$ and $C_u \in [0, 1]$. Parameters $C_E$ and $C_u$ are also obtained by calibration, while $C_H(t)$ is determined by employing the data of the infection transmissions within hospitals, usually available, as
                \begin{equation*}
                    C_H(t) = \frac{\alpha_H\left(\frac{\beta_I}{\gamma_I(t)} + \frac{\beta_E}{\gamma_E(t)} + (1 - \theta(t)) \frac{\beta_{I_u}(t)}{\gamma_{I_u}(t)}\right)}{(1 - \alpha_H)\beta_I\theta(t)\left(\left(1 - \frac{\omega(t)}{\theta(t)}\right)\frac{1}{\gamma_{H_R}(t)} + \frac{\omega(t)}{\theta(t)}\frac{1}{\gamma_{H_R}(t)}\right)}
                \end{equation*}
                with $\alpha_H$ being the percentage of cases (healthcare workers) infected by individuals in compartments $H_R$ or $H_D$.
                
            \item The compartment transition rates $\gamma_E, \gamma_I, \gamma_{I_u}, \gamma_{H_R}, \gamma_{H_D} \in (0, +\infty)$. They are constructed in accordance to the recent literature, based on the average duration (in days) of an individual in each infectious compartment, denoted by $d_E$, $d_I$, $d_{I_u}$, $d_{H_R}$ and $d_{H_D}$. Lets assume that $d_{I_u} = d_{H_R}$ and $d_{H_D} = d_{H_R} + \delta_R, \delta_R>0$. Thus,
                \begin{equation*}
                    \gamma_I = \frac{1}{d_E}, \quad \gamma_I(t) = \frac{1}{d_I - g(t)}, \quad \gamma_{I_u}(t) = \gamma_{H_R}(t) = \frac{1}{d_{I_u} + g(t)}, \quad \gamma_{H_D}(t) = \frac{1}{d_{I_u} + g(t) + \delta_R},
                \end{equation*}
                where $g(t) = d_g(1 - m_I(t))$ represents the decrease of the duration of $d_I$ due to the application of control measures at time $t$, and $d_g$ is the maximum number of days that $d_I$ can be decreased due to the control measures. The parameter $\delta_R$ is included in the calibration.
                
    \end{itemize}

        Our aim is to introduce stochasticity to the simplified $\theta$-SEIHRD model above. For this purpose, to keep the interpretability of the model, we will add some randomness to a set of parameters. In our approach, we will add randomness on the disease contact rates, $\beta$'s. More precisely, as in the simplified deterministic version proposed in \cite{ivorra2020} all the $\beta$'s depend on $\beta_I$, as indicated in \eqref{eq:betas}, it is sufficient to add uncertainty in the model parameter $\beta_I$ by turning it into a stochastic process, i.e a random variable at each time instant.
        
        Thus, we start by writing the model in terms of $\beta_I$. From the Equation \eqref{eq:betas}, we have the following relations:
        \begin{equation*}
            \beta_E = \beta_I A_E, \;\;\; \beta_{I_u} = \beta_I A_{I_u}, \;\;\; \beta_{H_R} = \beta_I A_{H_R}, \;\;\; \beta_{H_D} = \beta_I A_{H_D},
        \end{equation*}
        where
        \begin{equation*}
        \begin{aligned}
            A_E &= C_E, \\
            A_{I_u}(t) &= C_u + \frac{(1 - C_u)(1 - \theta(t))}{1 - \omega(t)}, \\
            A_{H_R}(t) &= A_{H_D}(t) = \frac{\alpha_H\left( \frac{1}{\gamma_I(t)} + \frac{A_E}{\gamma_E} + (1 - \theta(t))\frac{A_{I_u}(t)}{\gamma_{I_u(t)}} \right)}{(1 - \alpha_H)\theta(t)\left( (1 - \frac{\omega(t)}{\theta(t)})\frac{1}{\gamma_{H_R}(t)} + \frac{\omega(t)}{\theta(t)}\frac{1}{\gamma_{H_D}(t)} \right)}.
        \end{aligned}
        \end{equation*}
        
        By using the previous notation, the simplified $\theta$-SEIHRD model in \eqref{eq:thetaSEIHRD_model} can be rewritten as,
        \begin{equation}\label{eq:theta_seihrd_2}
        \begin{aligned}
            \frac{\d S}{\d t}(t) &= -\beta_I\frac{S(t)M(t)}{N}, \\
            \frac{\d E}{\d t}(t) &= \beta_I\frac{S(t)M(t)}{N} - \gamma_E E(t), \\
            \frac{\d I}{\d t}(t) &= \gamma_E E(t) - \gamma_I I(t), \\
            \frac{\d I_u}{\d t}(t) &= (1 - \theta(t))\gamma_I I(t) - \gamma_{I_u} I_u(t), \\
            \frac{\d H_R}{\d t}(t) &= \theta(t)\left(1 - \frac{\omega(t)}{\theta(t)}\right)\gamma_I I(t) - \gamma_{H_R} H_R(t), \\
            \frac{\d H_D}{\d t}(t) &= \omega(t)\gamma_I I(t) - \gamma_{H_D} H_D(t), \\
            \frac{\d R_d}{\d t}(t) &= \gamma_{H_R}(t) H_R(t), \\
            \frac{\d R_u}{\d t}(t) &= \gamma_{I_u}(t) I_u(t), \\
            \frac{\d D}{\d t}(t) &= \gamma_{H_D}(t) H_D(t), 
        \end{aligned}
        \end{equation}
        where
        \begin{equation}\label{eq:M}
            M(t) = m_E A_E E(t) + m_I I(t) + m_{I_u} A_{I_u} I_u(t) + m_{H_R} A_{H_R} H_R(t) + m_{H_D} A_{H_D} H_D(t).
        \end{equation}
        
        We now incorporate a stochastic component into the model by, as mentioned, by replacing the constant parameter $\beta_I$ of the deterministic model by a stochastic process.
        Lets then assume that, instead of considering a constant value $\beta_I$, the disease contact rate in compartment $I$ follows a newly introduced stochastic process $\tilde{\beta}_I(t)$. In order to preserve the positiveness in the parameter definition (imposed by its biological interpretation), we choose the well-known CIR process \cite{cox1985}. The main advantage of the CIR process is that it theoretically ensures the spacial states to be non-negative. Further, the CIR process is a mean-reverting process meaning that, in the long term, the dynamics of the process tend to a prescribed value (the average). This can have a biological interpretation since, at the early stages of the disease, the contacts between people are less controlled, so more volatile. In the mid-long term, the individual contacts and the disease spread patterns are more studied and the variability in its estimation ranges is reduced. Accordingly, the dynamics of $\tilde{\beta}_I$ read
            \begin{equation}\label{eq:betaI_CIR}
                \d \tilde{\beta}_I(t) = \nu_{\beta_I}(\mu_{\beta_I} - \tilde{\beta}_I(t))\dt + \sigma_{\beta_I}\sqrt{\tilde{\beta}(t)}\dW(t)
            \end{equation}
            where $\nu_{\beta_I}$ is the mean reverting speed, $\mu_{\beta_I}$ is the long-term average, $\sigma_{\beta_I}$ is the volatility and $\dW(t)$ is a Brownian motion increment.
        
        Therefore, the corresponding system of SDEs governing the stochastic $\theta$-SEIHRD model is given by,
        \begin{equation}\label{eq:sto_theta_seihrd}
        \begin{aligned}
            \d S(t) &= \tilde{\beta}_I(t)\frac{S(t)M(t)}{N}\dt \\
            \d E(t) &= \left(\tilde{\beta}_I(t)\frac{S(t)M(t)}{N} - \gamma_E E(t)\right)\dt \\
            \d I(t) &= \left(\gamma_E E(t) - \gamma_I I(t)\right)\dt, \\
            \d I_u(t) &= \left((1 - \theta(t))\gamma_I I(t) - \gamma_{I_u} I_u(t)\right)\dt, \\
            \d H_R(t) &= \left(\theta(t)\left(1 - \frac{\omega(t)}{\theta(t)}\right)\gamma_I I(t) - \gamma_{H_R} H_R(t)\right)\dt, \\
            \d H_D(t) &= \left(\omega(t)\gamma_I I(t) - \gamma_{H_D} H_D(t)\right)\dt, \\
            \d R_d(t) &= \gamma_{H_R}(t) H_R(t) \dt, \\
            \d R_u(t) &= \gamma_{I_u}(t) I_u(t) \dt, \\
            \d D(t) &= \gamma_{H_D}(t) H_D(t) \dt, \\
            \d \tilde{\beta}_I(t) &= \nu_{\beta_I}(\mu_{\beta_I} - \tilde{\beta}_I(t))\dt + \sigma_{\beta_I}\sqrt{\tilde{\beta}(t)}\dW(t).
        \end{aligned}
        \end{equation}
        with $M(t)$ as defined in Equation \eqref{eq:M}. Although only one source of randomness is introduced, the solution of the whole system becomes a set of stochastic processes, due to the dependence of the remaining equations on the first equations for $S$ and $E$, and their own dependence on $\tilde{\beta}_I$. Note also that we could add more time dependency to $\tilde{\beta}_I$ both in the deterministic model (which is not the case in \cite{ivorra2020}) as well as in the stochastic version (by considering either $\nu_{\beta_I}$, $\mu_{\beta_I}$ or  $\sigma_{\beta_I}$ time dependent).
    
        Note that our approach could be generalized to a setting where some of the parameters were independent, in this case we would consider each one as a Gaussian random variable and possible correlations between some or all of them. In this more general setting, a certain number of different (possibly correlated) Brownian motion processes would come into place.
        
    \section{Numerical solution}
    
        For a set of constant initial data $S(0),\, E(0), \, I(0), \, I_u(0), \, H_R(0), \, H_D(0),\, R_d(0),\, R_u(0)$, $D(0)$ and $\tilde{\beta}_I(0)$, the system of SDEs (\ref{eq:sto_theta_seihrd}) has a unique strong solution. This follows from the fact that the coefficients of the SDEs are locally Lipchitz continuous and the initial condition is a constant value (see \cite{kloeden1992}, for example).
        
        Although in the stochastic $\theta$-SEIHRD model (\ref{eq:sto_theta_seihrd}) the solution is a finite set of stochastic processes, we maintain the same notation we used for the finite set of real valued functions representing the solution of the deterministic $\theta$-SEIHRD model.
    
        As the system of SDEs (\ref{eq:sto_theta_seihrd}) is nonlinear, it is not possible to obtain a closed-form expression for the solution, as it also happens with the deterministic version. Therefore, the use of numerical methods for solving (\ref{eq:sto_theta_seihrd}) becomes mandatory. 
        Here, we adopt the following strategy. First, we perform a simulation of the dynamics of $\tilde{\beta}_I(t)$, in accordance with the CIR process. After, we solve the resulting ODE system for each simulated path of $\tilde{\beta}_I(t)$. In this way, we obtain a set of random walks for each stochastic process representing a model variable.
        
        The CIR process is a well-studied dynamics often employed in computational finance (see \cite{oosterlee2019} for example), which satisfies the SDE in Equation \eqref{eq:betaI_CIR}. From the mathematical point of view, one of its relevant features is that the underlying distribution is known analytically, relying on the \emph{non-central chi-squared} distribution. Thus, given two time points, $s$ and $t$, $s<t$, the conditional distribution of $\tilde{\beta}_I$ defined in Equation \eqref{eq:betaI_CIR} reads
            \begin{equation*}
                \tilde{\beta}_I(t)|\tilde{\beta}_I(s) \sim c(t,s)\cdot\chi^2\left(d, \frac{\e^{-\nu_{\beta_I}(t - s)}}{c(t, s)}\tilde{\beta}_I(s) \right),
            \end{equation*}
            where
            \begin{equation*}
                c(t, s) = \frac{\sigma_{\beta_I}^2}{4\mu_{\beta_I}}\left(1 - \e^{-\nu_{\beta_I}(t - s)}\right), \quad d = \frac{4\nu_{\beta_I}\mu_{\beta_I}}{\sigma_{\beta_I}^2},
            \end{equation*}
            and
            $\chi^2(a, b)$ is the non-central chi-squared distribution with $a$ degrees of freedom and non-centrality parameter $b$. This form the basis for an \emph{exact simulation} scheme, which can be used to obtain realizations of $\tilde{\beta}_I$. Given a set of $m+1$ time points, $\{t_i\}_0^m$, where the solution will be computed, we have, for $i=0,\dots,m-1$,
            \begin{equation*}
            \begin{aligned}
                c(t_{i+1}, t_i) &= \frac{\sigma_{\beta_I}^2}{4\mu_{\beta_I}}\left(1 - \e^{-\nu_{\beta_I}(t_{i+1} - t_i)}\right), \\
                \tilde{\beta}_I(t_{i+1}) &= c(t_{i+1}, t_i)\chi^2\left(d, \frac{\e^{-\nu_{\beta_I}(t_{i+1} - t_i)}}{c(t_{i+1} - t_i)}\tilde{\beta}_I(t_i) \right)
            \end{aligned}
            \end{equation*}
            with the constant parameter $d = \frac{4\nu_{\beta_I}\mu_{\beta_I}}{\sigma_{\beta_I}^2}$ and some initial value $\tilde{\beta}_I(t_0) = \tilde{\beta}_I(0)$. By employing this scheme, we generate $n$ simulated discrete sample paths of $\tilde{\beta}_I$.
        
        Once the sample paths of $\tilde{\beta}_I$ have been obtained, we numerically solve $n$ ODE systems, one for each of these paths. For this purpose, we choose the explicit Runge-Kutta method of order 5(4), known as \emph{RK45}, \emph{RKDP} or Dormand–Prince method, see \cite{dormand1980}. Note that the time discretization needs to include the time points used in the simulation of $\tilde{\beta}_I$. Alternatively, the whole system in Equation \eqref{eq:sto_theta_seihrd} can be solved by the Euler-Maruyama scheme (see \cite{kloeden1992}, for example), which is the stochastic version of the classical Euler scheme for solving ODE systems. Although this choice would provide faster solutions, it involves discretization errors, giving rise to some lack of precision and robustness.
        
\begin{figure}[h!]
            \centering
                \subfigure[$S(t)$]{\includegraphics[width=0.49\textwidth]{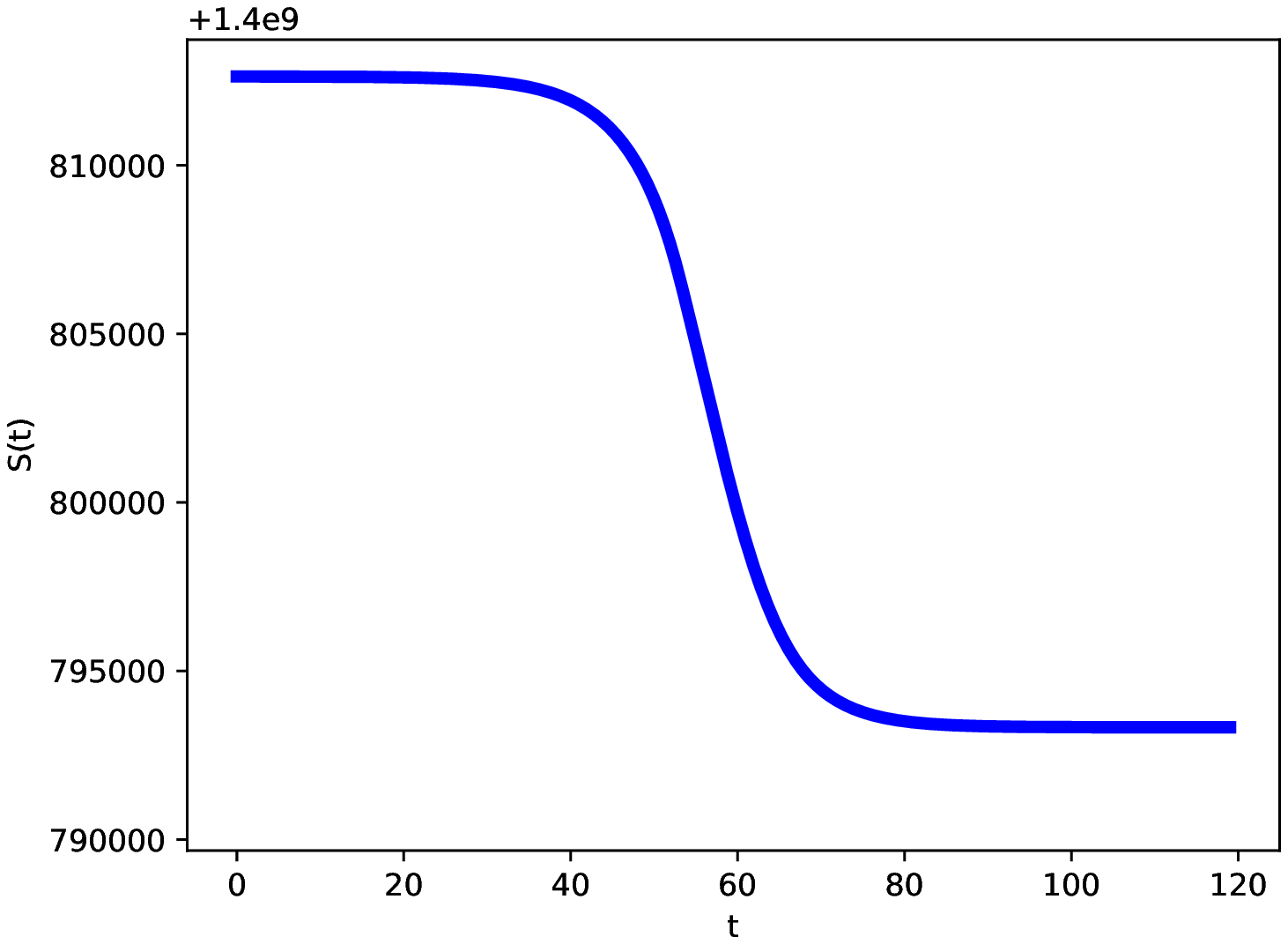}}
                \subfigure[$S(t)$]{\includegraphics[width=0.49\textwidth]{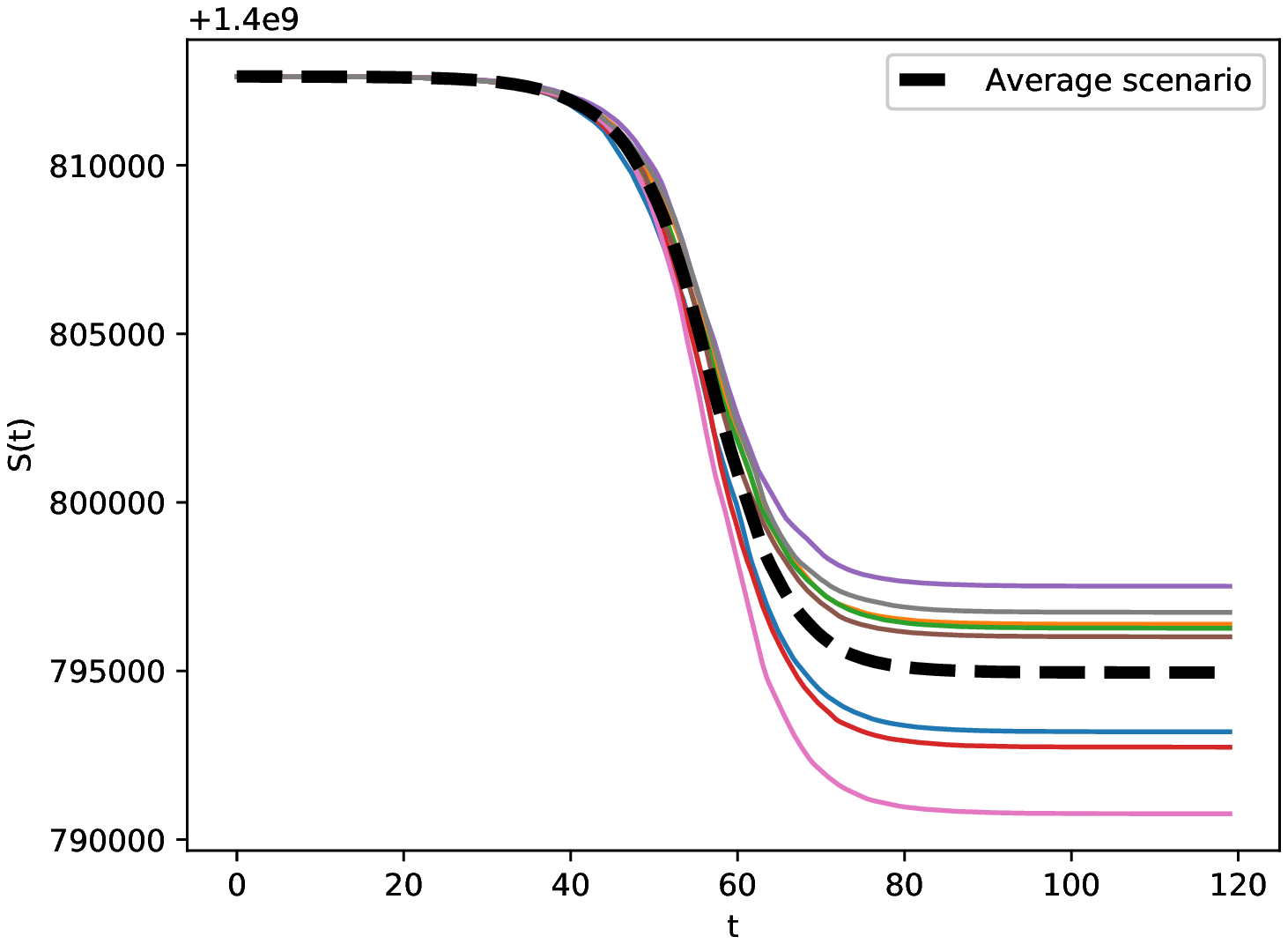}}
                \subfigure[$E(t)$]{\includegraphics[width=0.49\textwidth]{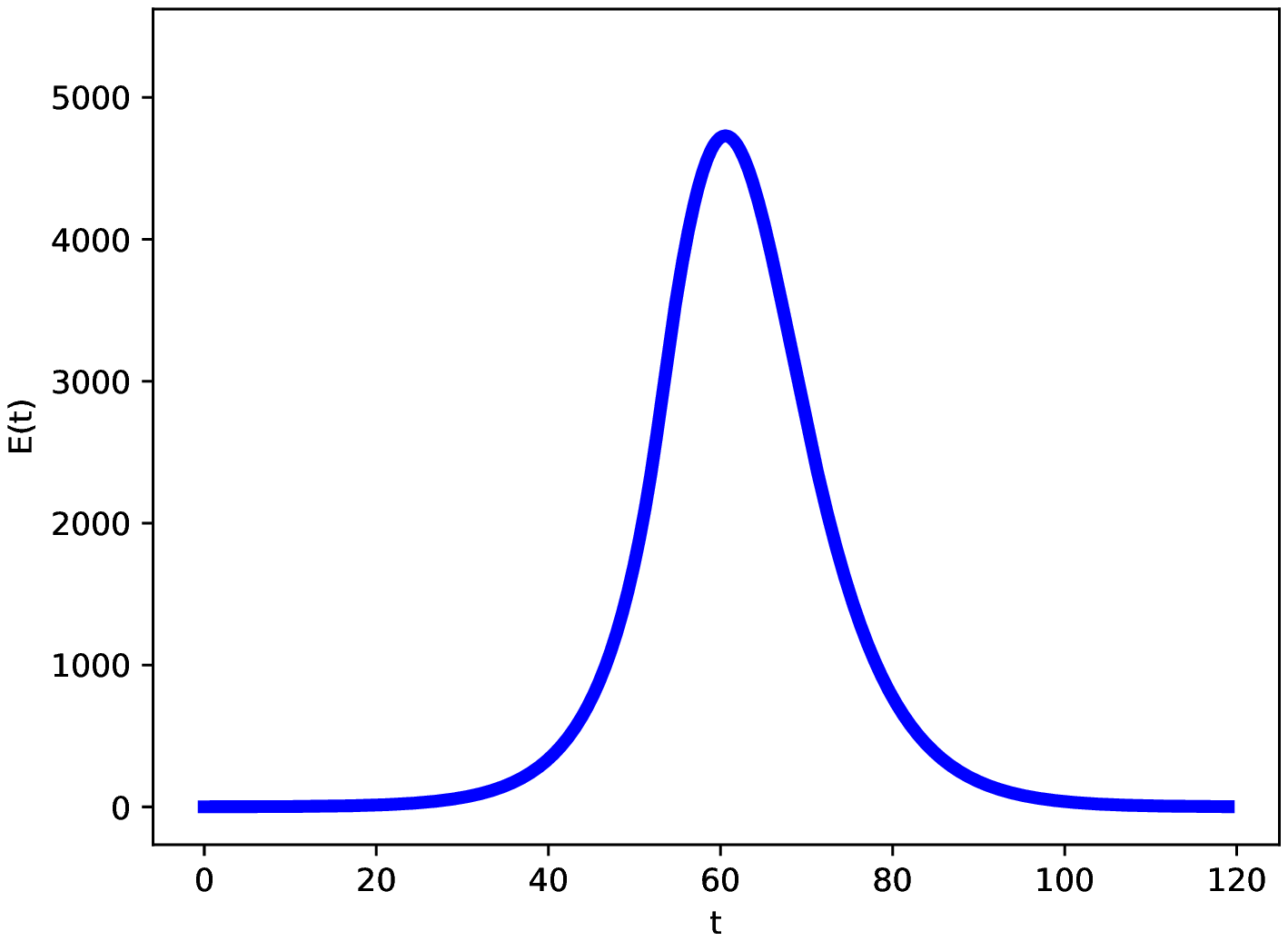}}
                \subfigure[$E(t)$]{\includegraphics[width=0.49\textwidth]{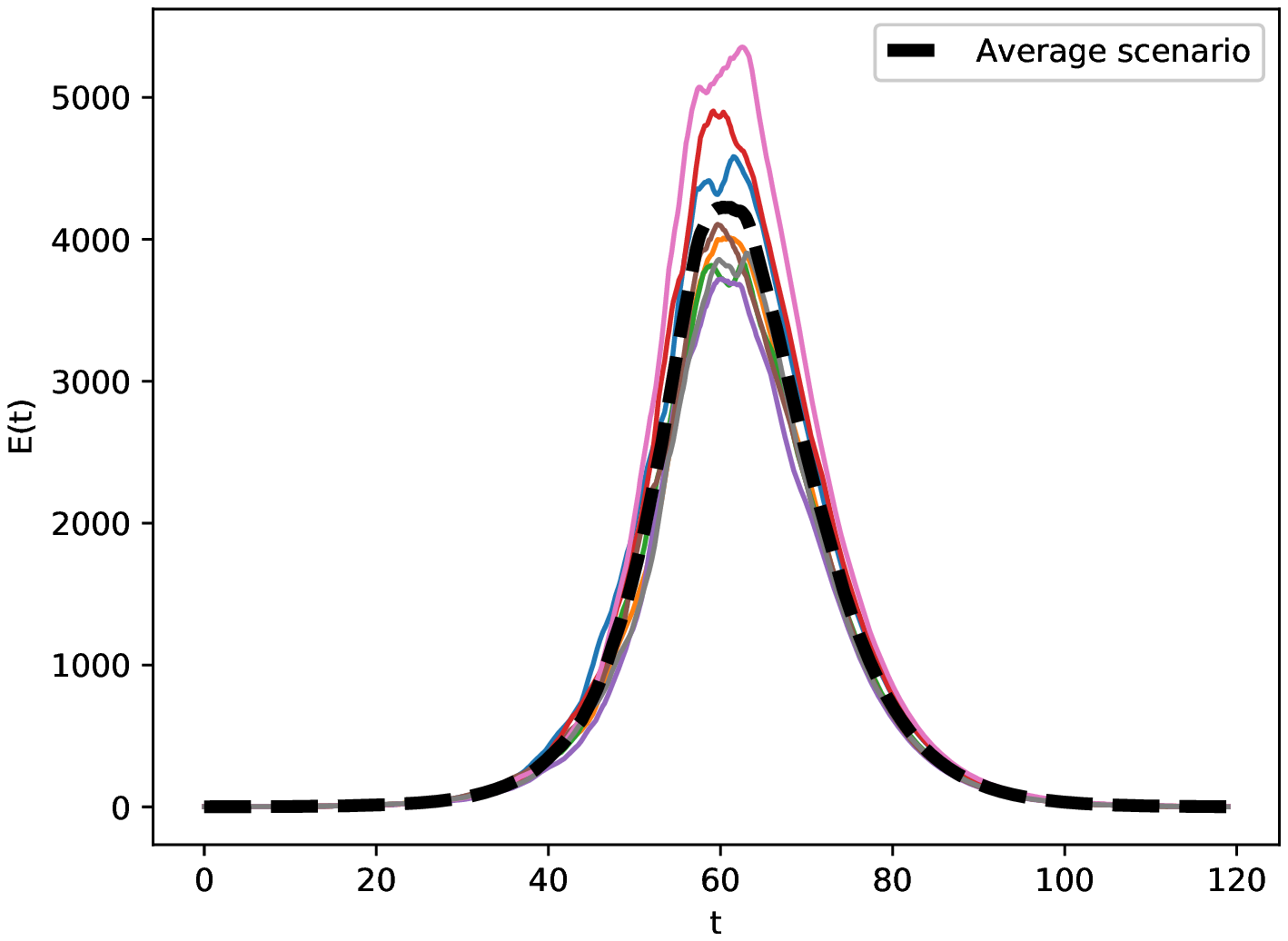}}
                \subfigure[$I(t)$]{\includegraphics[width=0.49\textwidth]{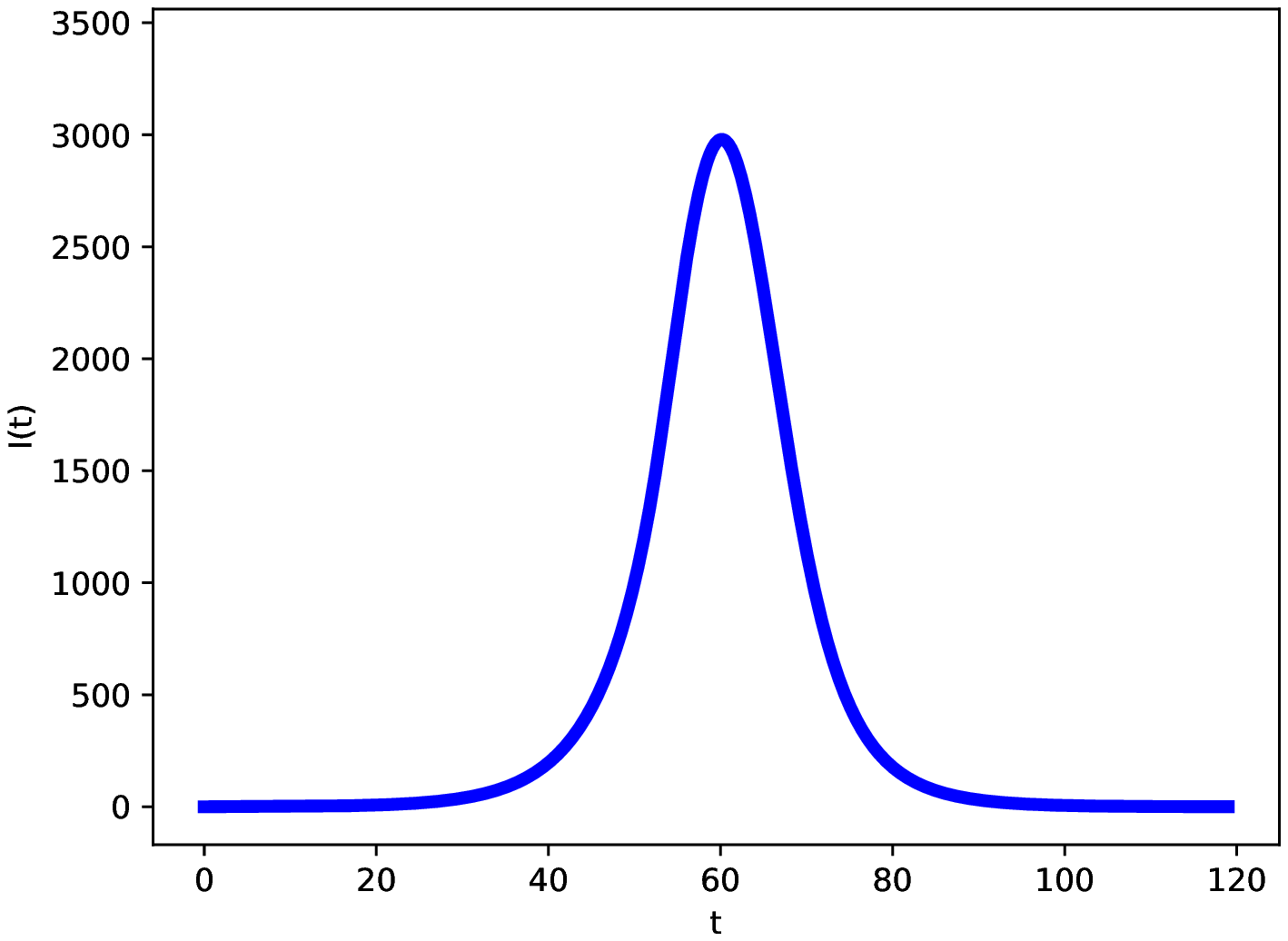}}
                \subfigure[$I(t)$]{\includegraphics[width=0.49\textwidth]{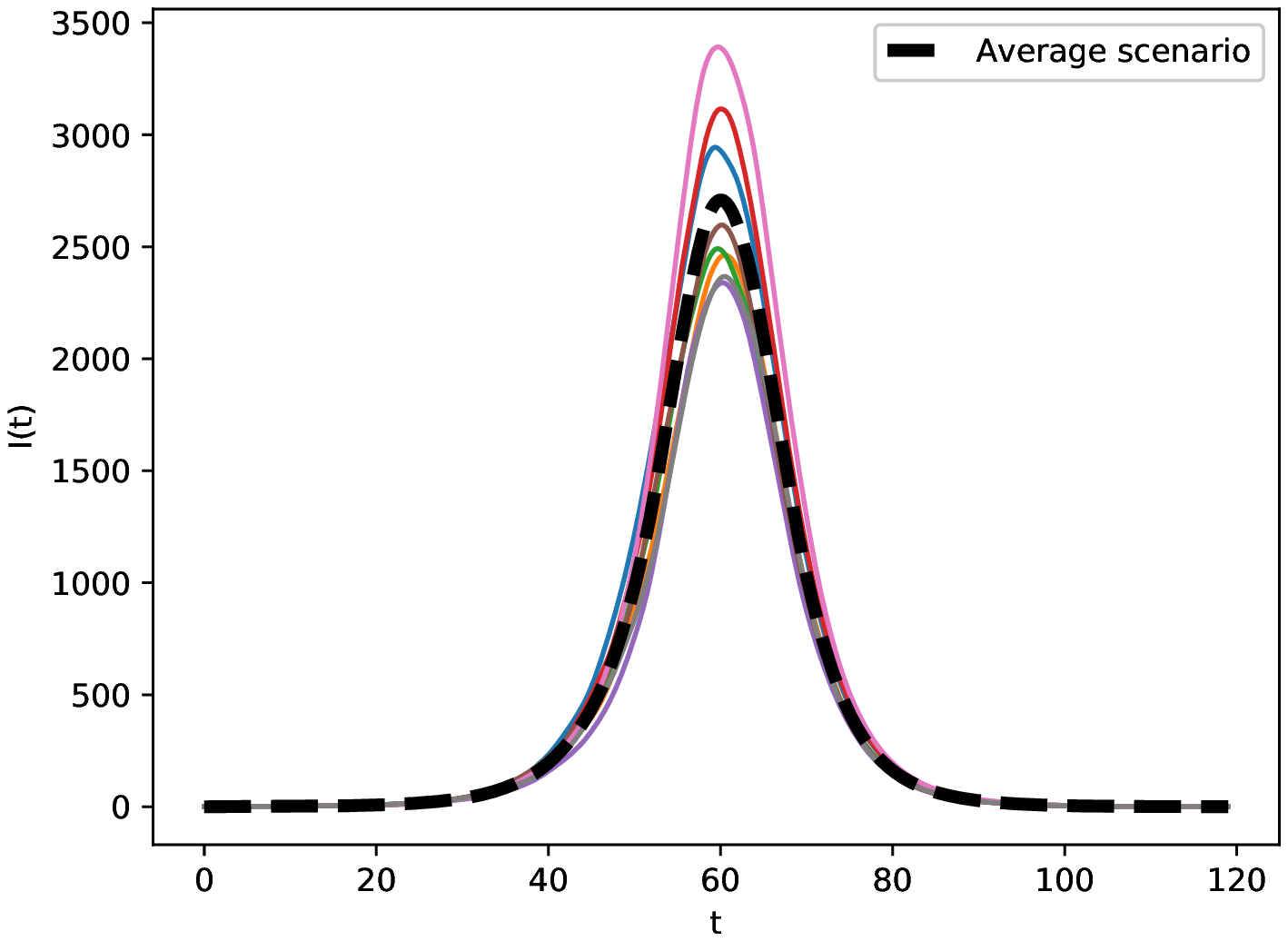}}
            \caption{Deterministic vs. Stochastic: $\nu_{\beta_I} = 1$, $\mu_{\beta_I} = \beta_I$ and $\sigma_{\beta_I} = 0.1$, with $n = 8$ Monte Carlo simulations.}
            \label{fig:Deterministic_vs_Stochastic_1}
        \end{figure}

\begin{figure}[h!]
            \centering
                \subfigure[$I_u(t)$]{\includegraphics[width=0.49\textwidth]{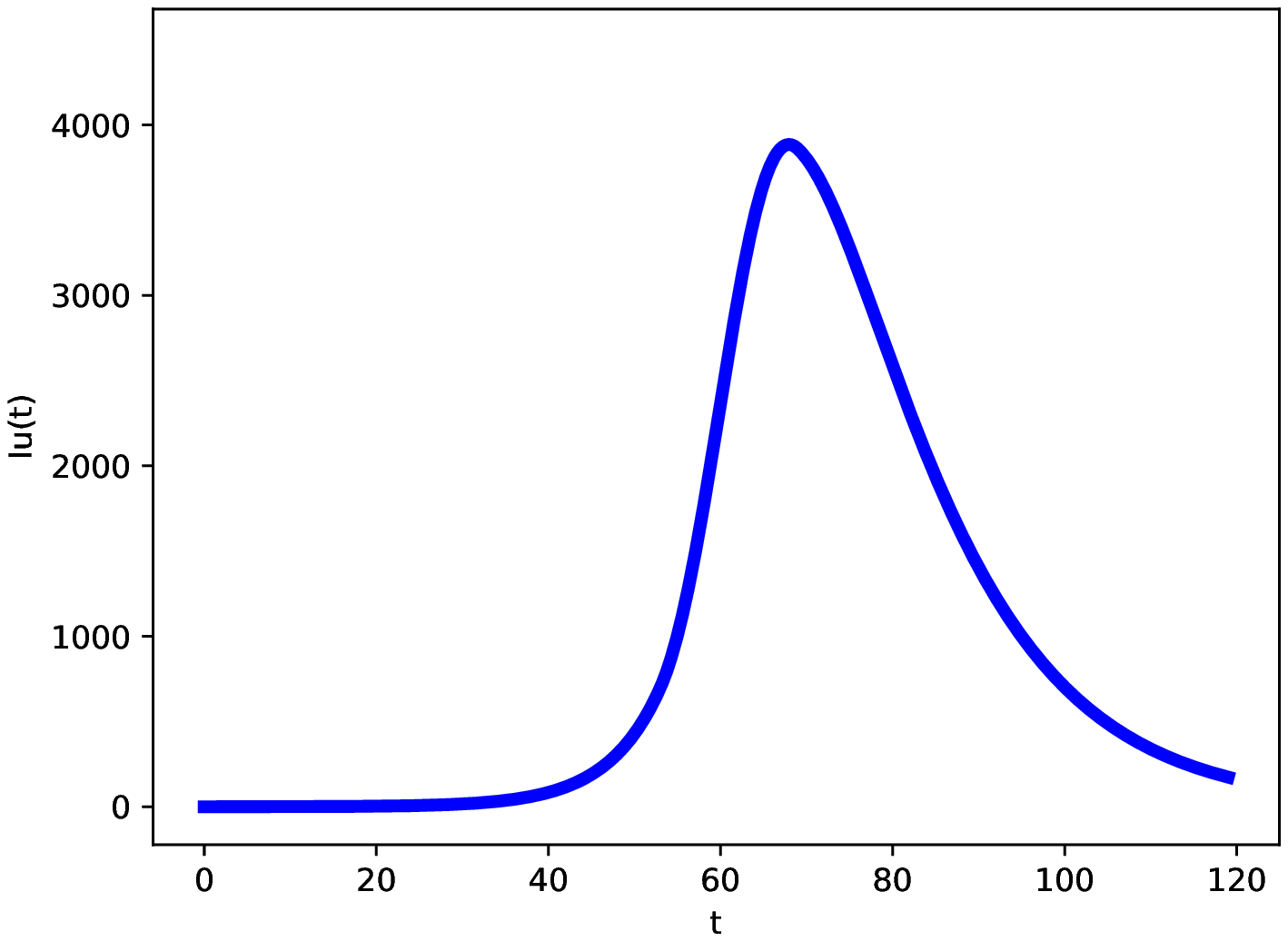}}
                \subfigure[$I_u(t)$]{\includegraphics[width=0.49\textwidth]{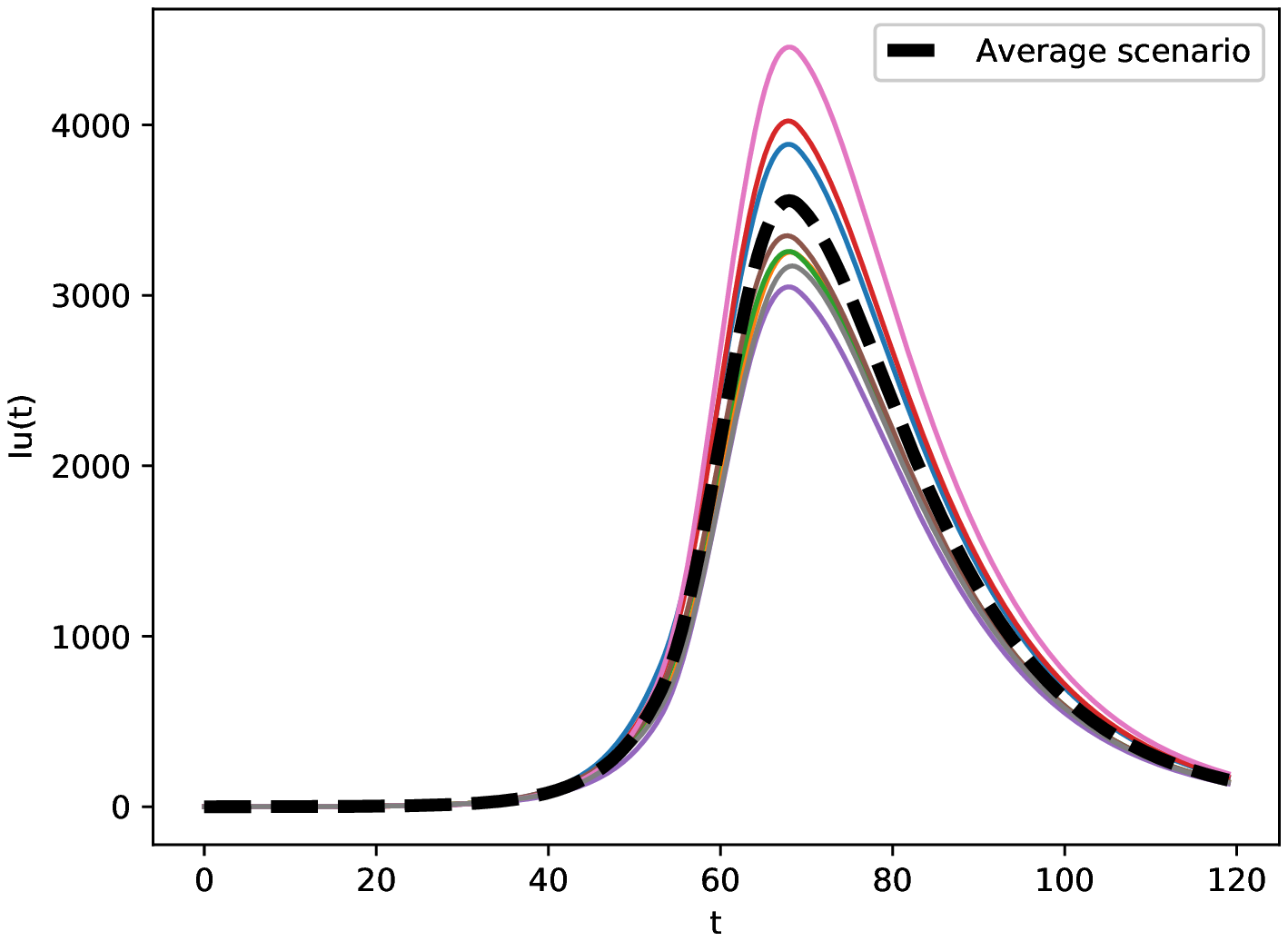}}
                \subfigure[$H_R(t)$]{\includegraphics[width=0.49\textwidth]{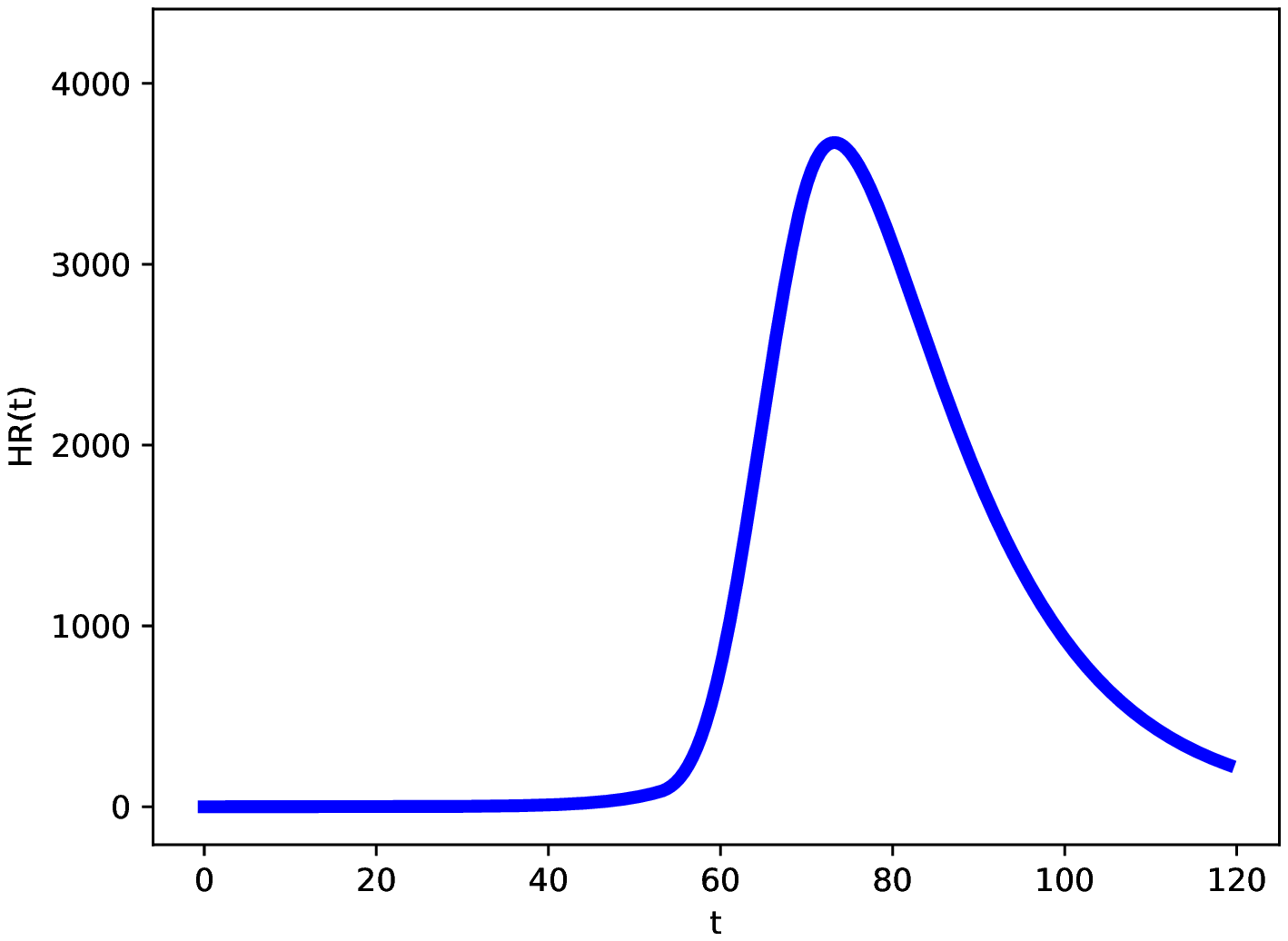}}
                \subfigure[$H_R(t)$]{\includegraphics[width=0.49\textwidth]{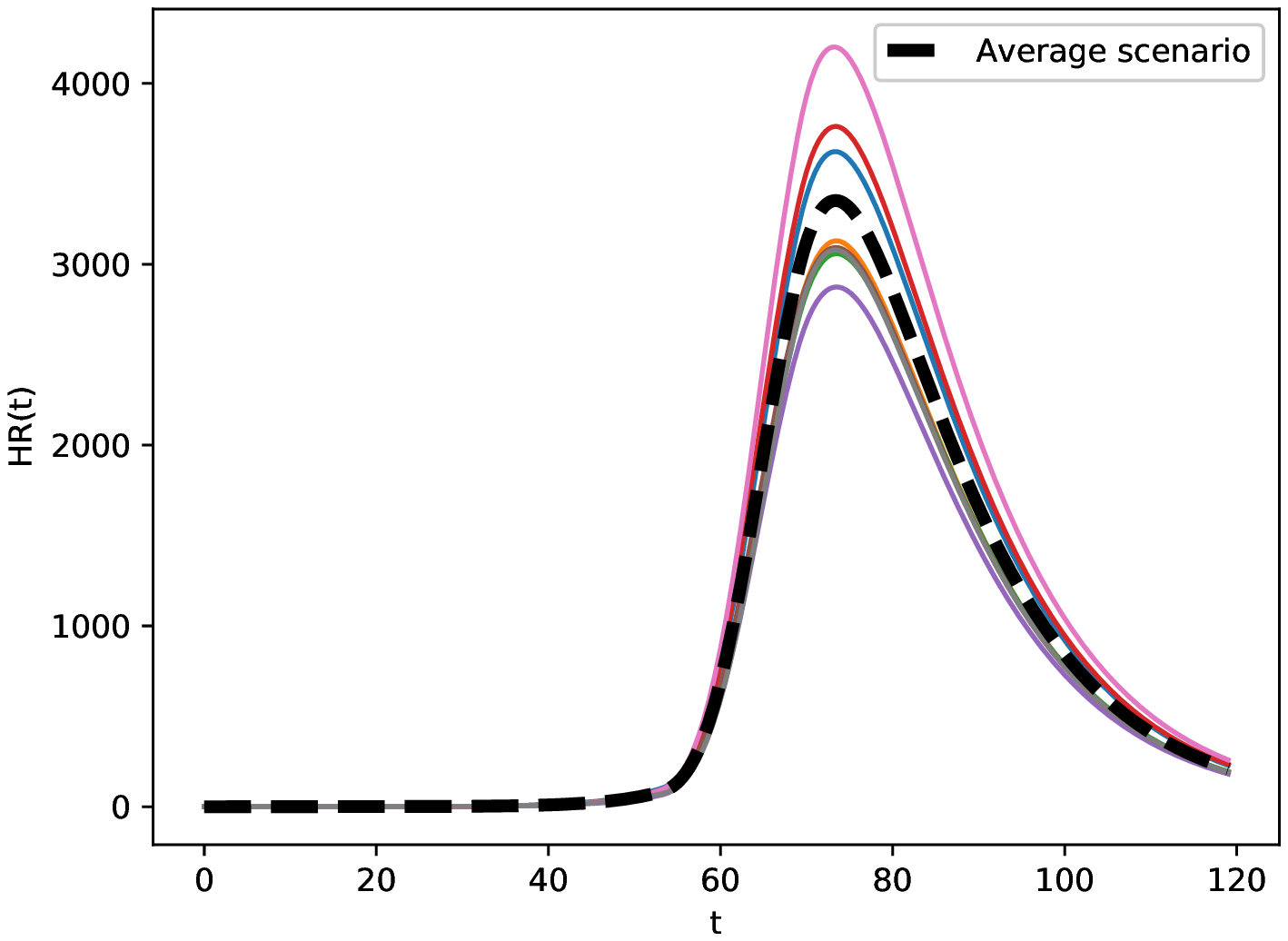}}
                \subfigure[$H_D(t)$]{\includegraphics[width=0.49\textwidth]{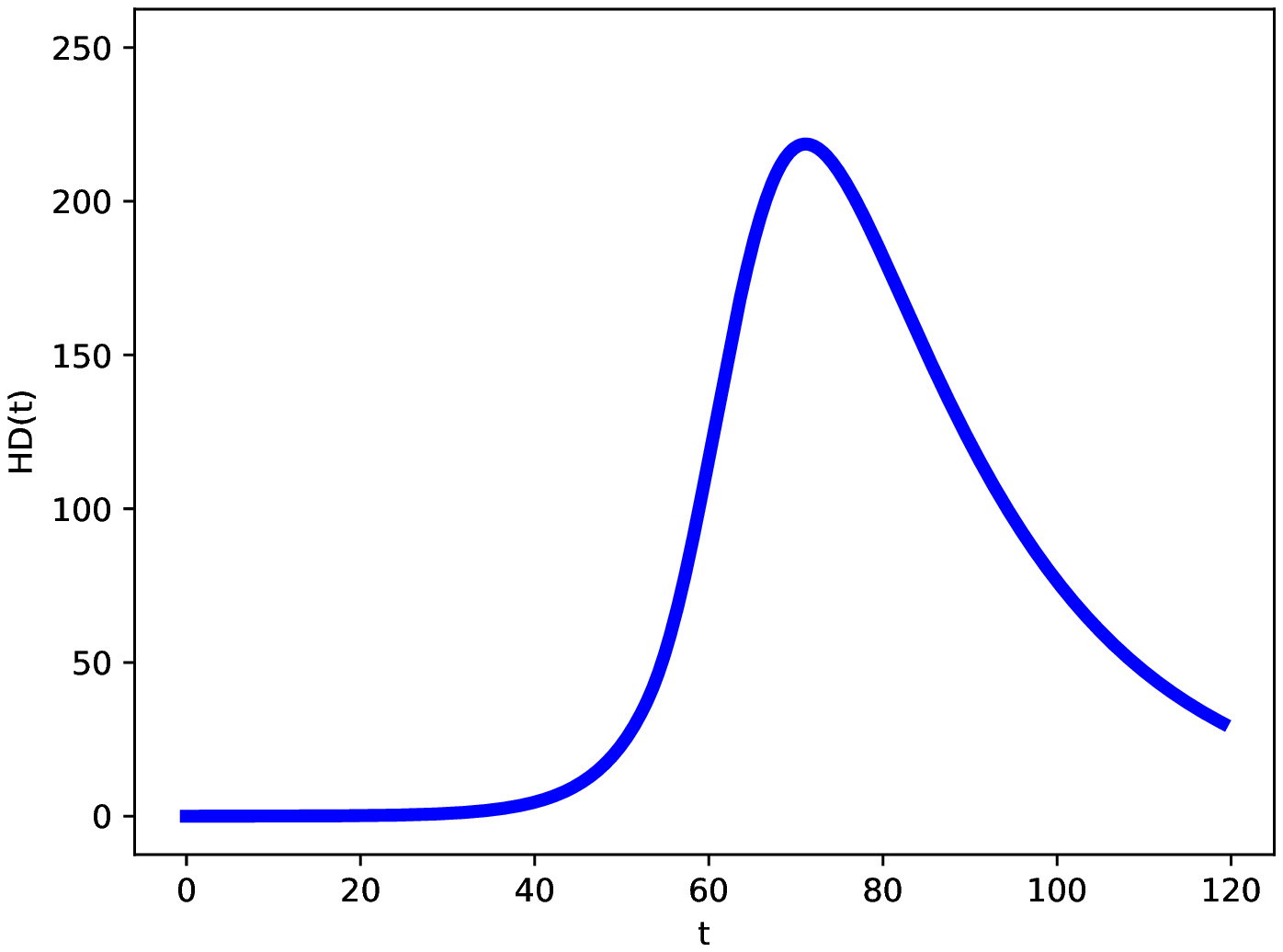}}
                \subfigure[$H_D(t)$]{\includegraphics[width=0.49\textwidth]{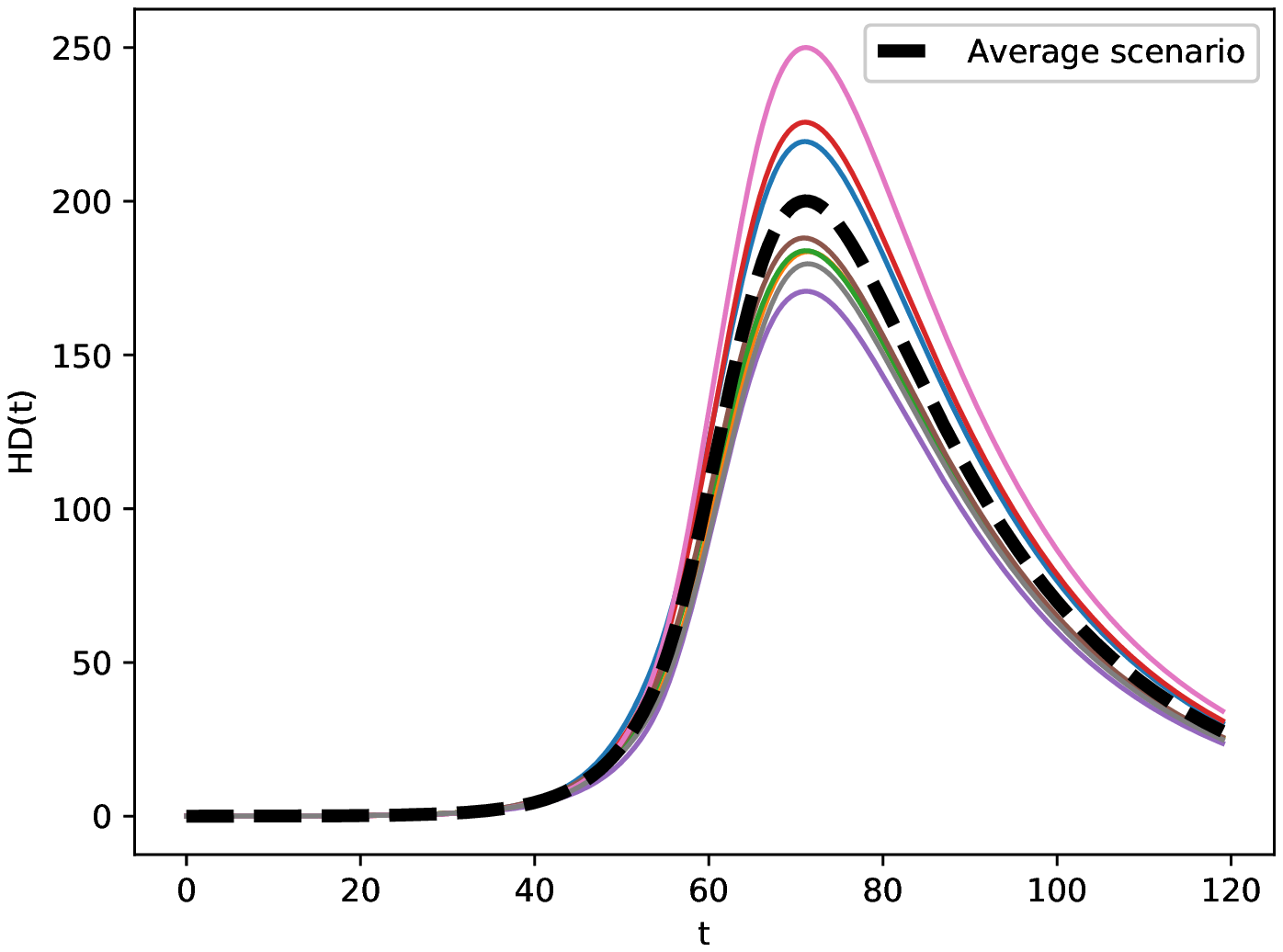}}
            \caption{Deterministic vs. Stochastic: $\nu_{\beta_I} = 1$, $\mu_{\beta_I} = \beta_I$ and $\sigma_{\beta_I} = 0.1$, with $n = 8$ Monte Carlo simulations.}
            \label{fig:Deterministic_vs_Stochastic_2}
        \end{figure}

\begin{figure}[h!]
            \centering
                \subfigure[$R_d(t)$]{\includegraphics[width=0.49\textwidth]{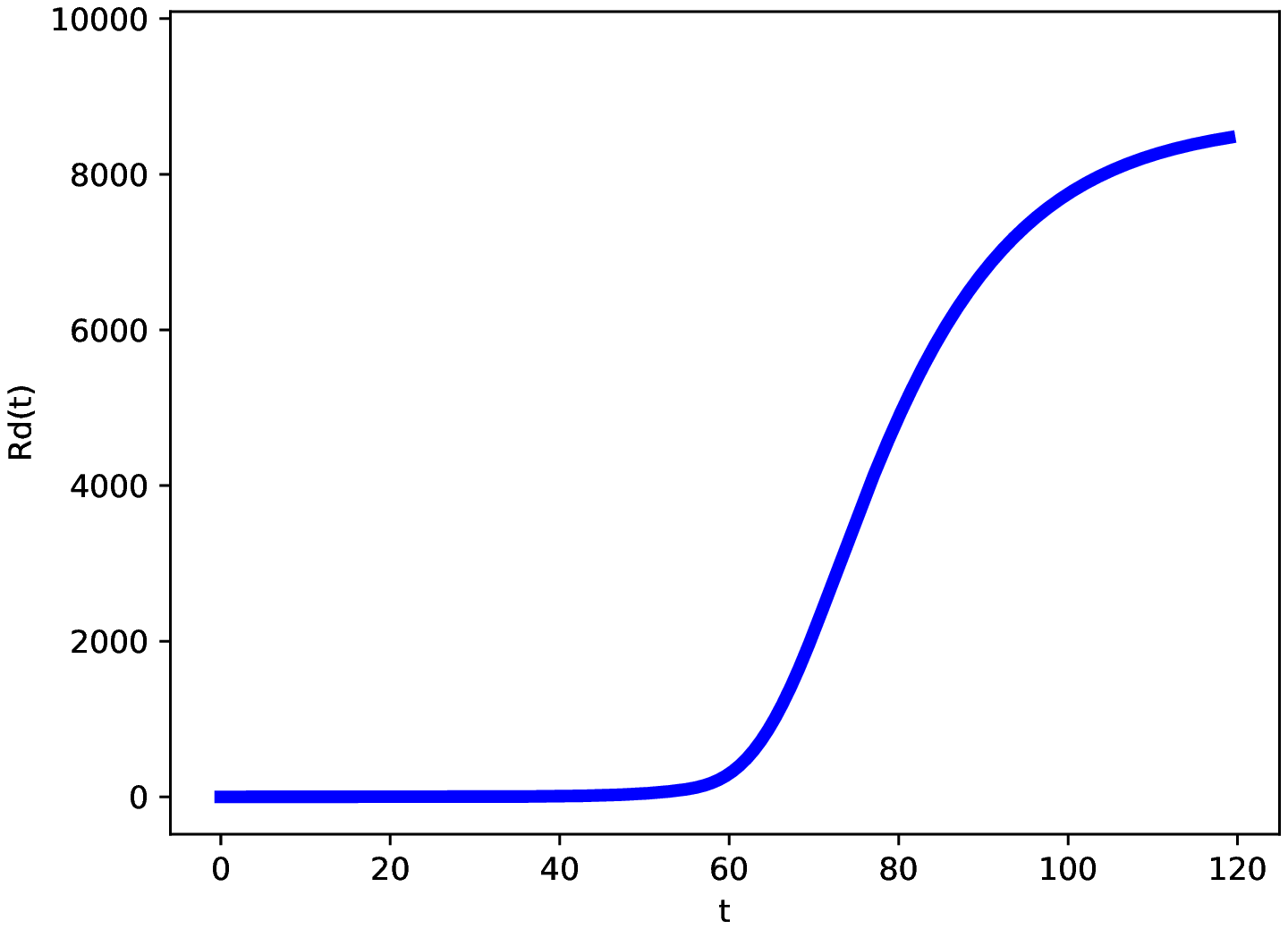}}
                \subfigure[$R_d(t)$]{\includegraphics[width=0.49\textwidth]{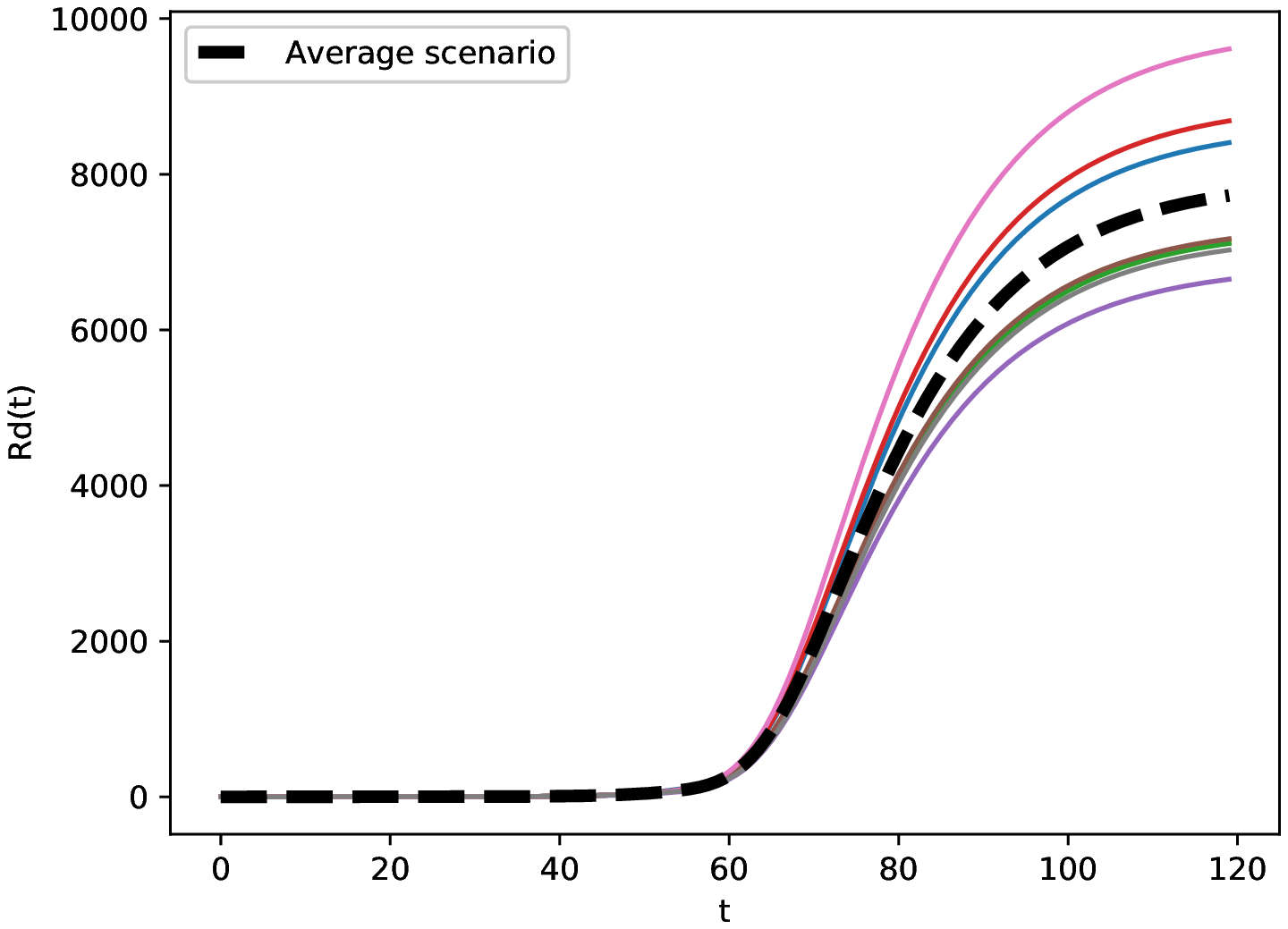}}
                \subfigure[$R_u(t)$]{\includegraphics[width=0.49\textwidth]{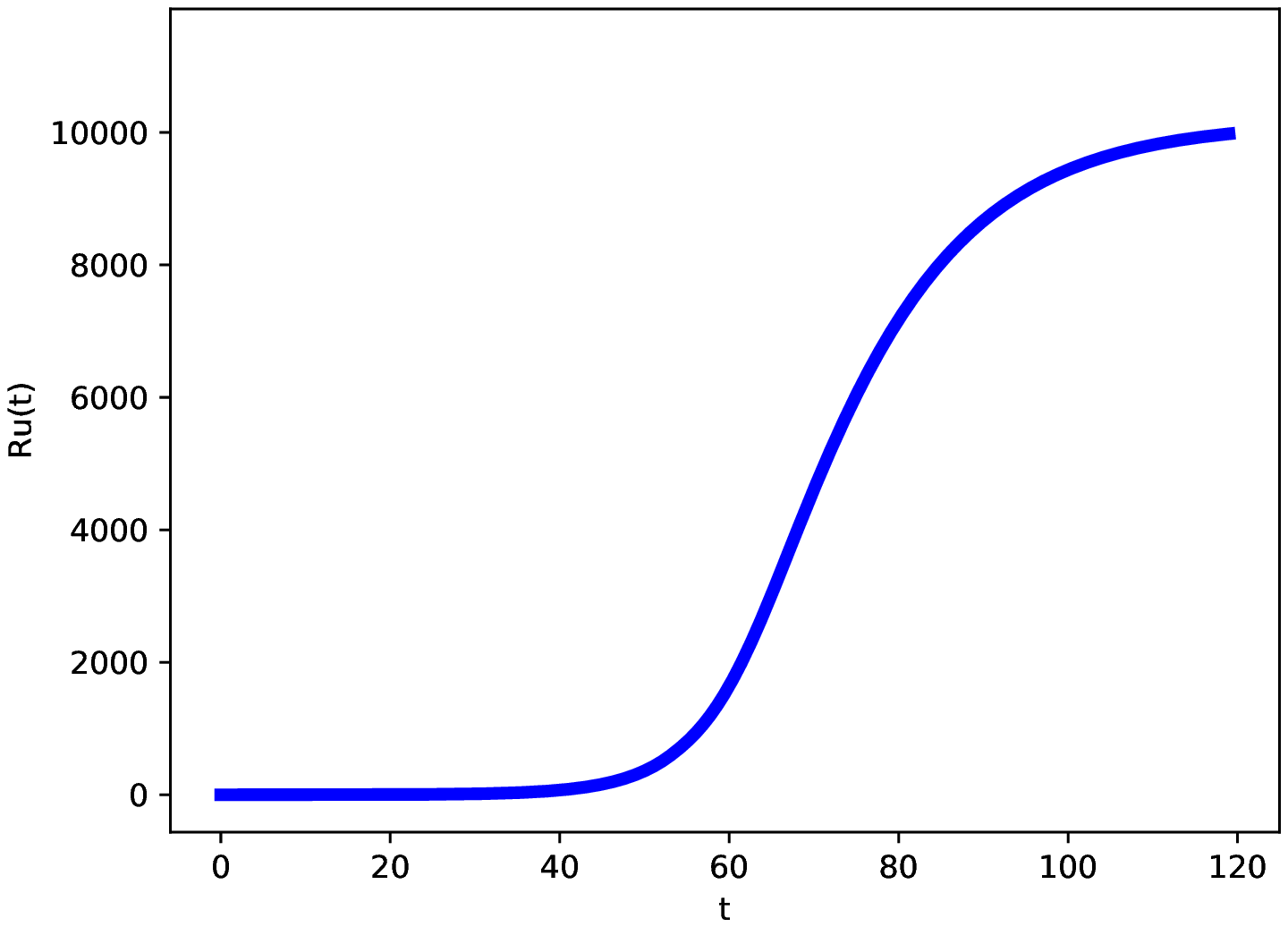}}
                \subfigure[$R_u(t)$]{\includegraphics[width=0.49\textwidth]{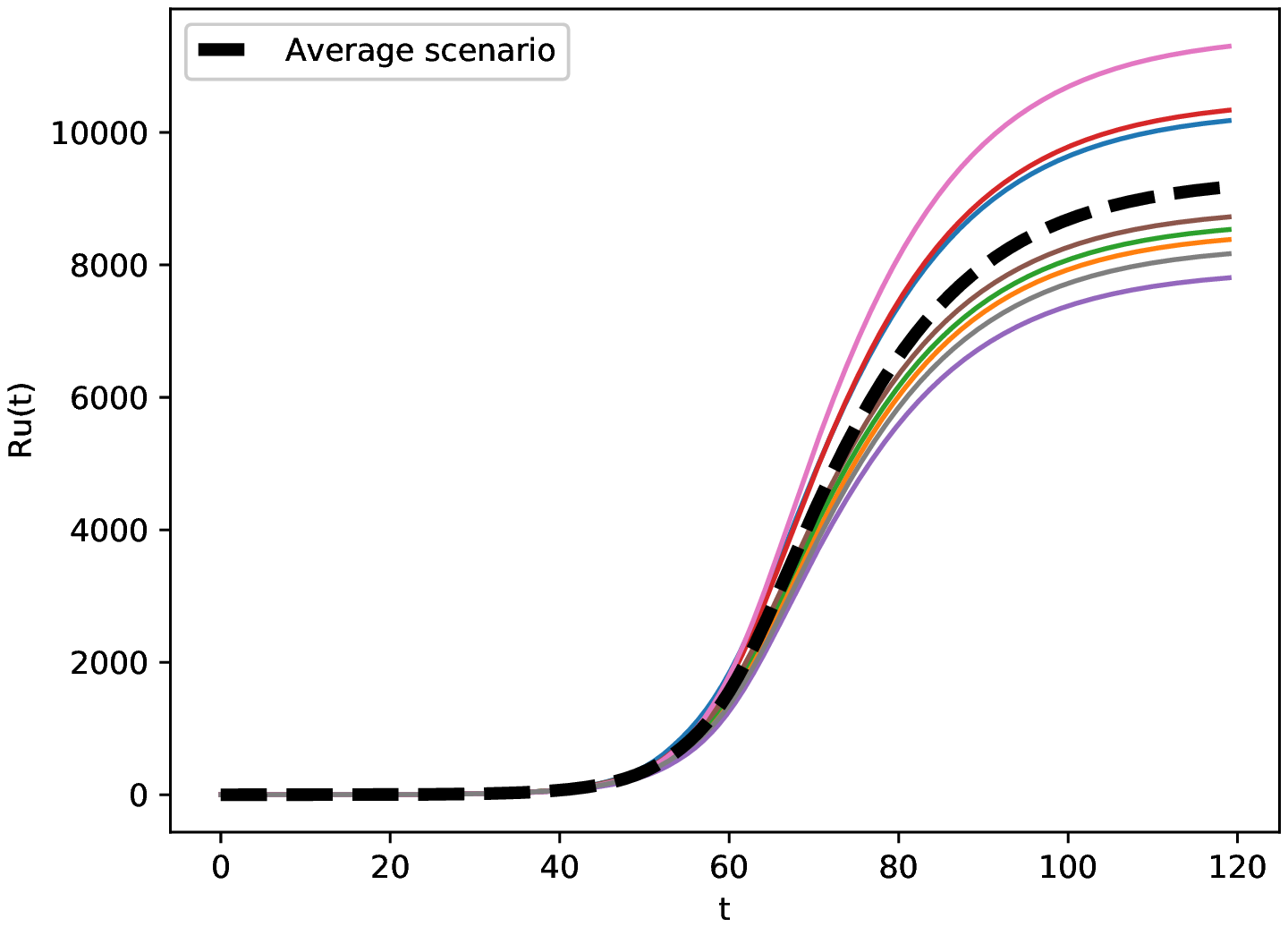}}
                \subfigure[$D(t)$]{\includegraphics[width=0.49\textwidth]{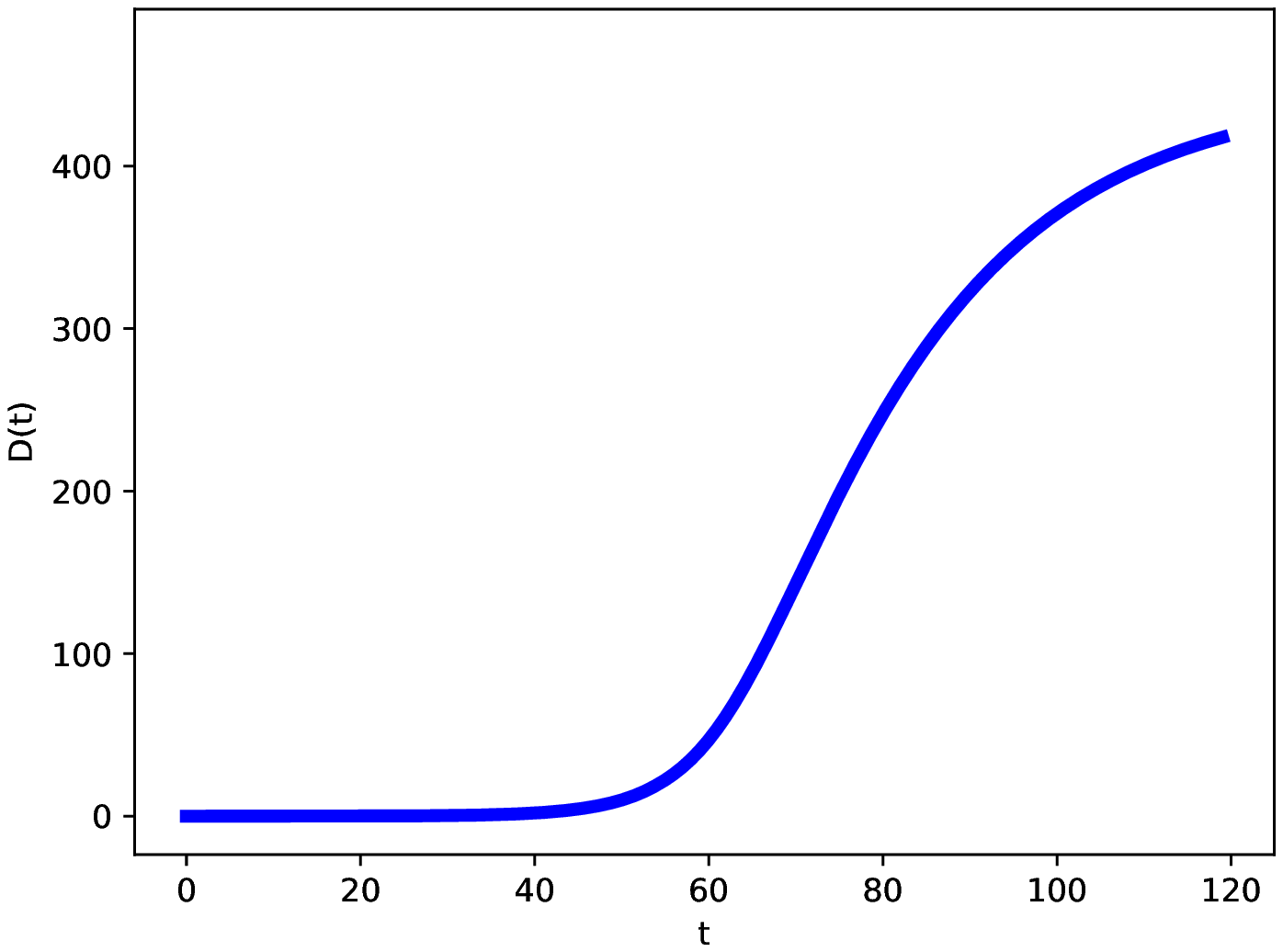}}
                \subfigure[$D(t)$]{\includegraphics[width=0.49\textwidth]{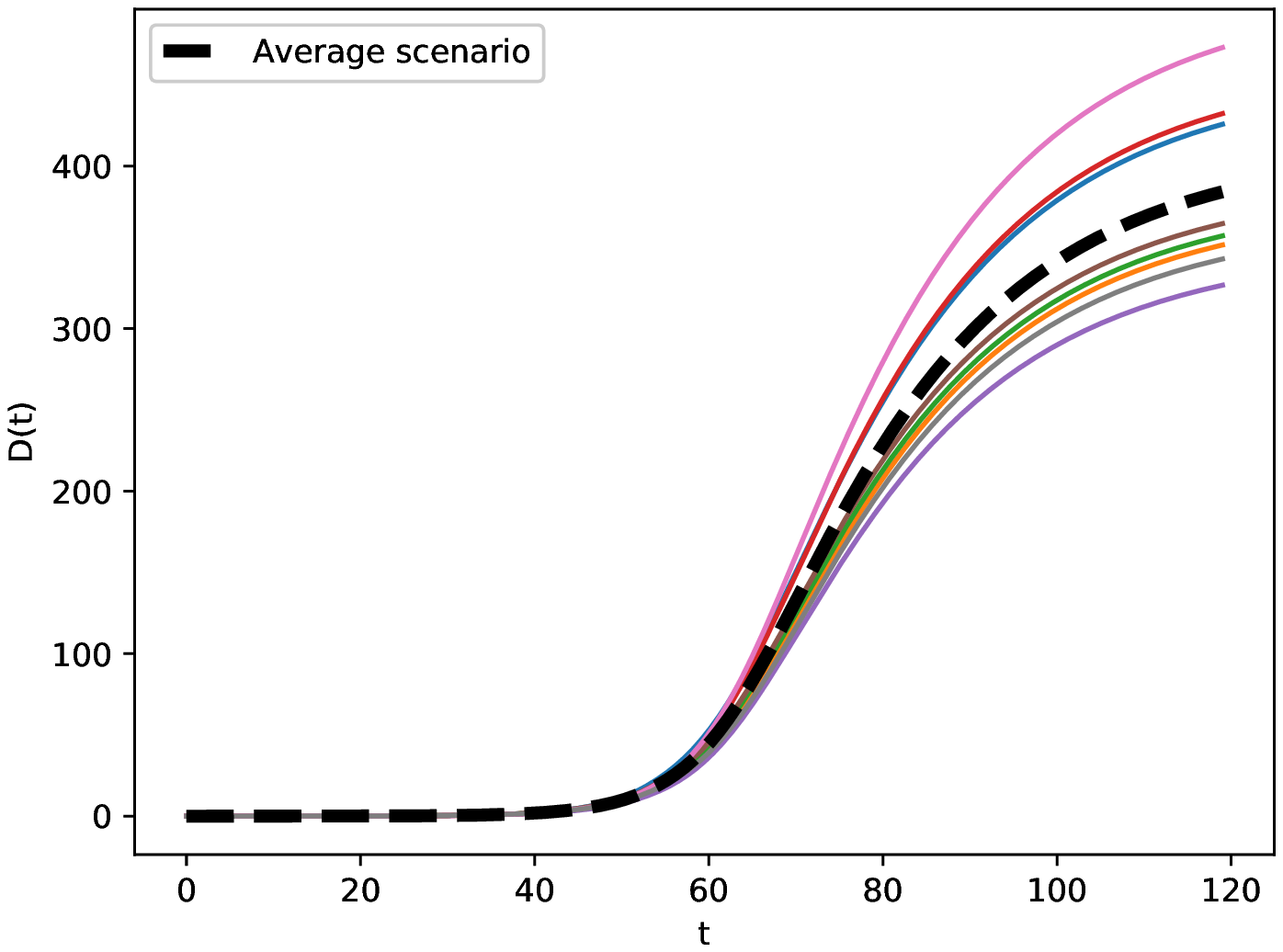}}
            \caption{Deterministic vs. Stochastic: $\nu_{\beta_I} = 1$, $\mu_{\beta_I} = \beta_I$ and $\sigma_{\beta_I} = 0.1$, with $n = 8$ Monte Carlo simulations.}
            \label{fig:Deterministic_vs_Stochastic_3}
        \end{figure}

        In order to illustrate the potential of the stochastic version of $\theta$-SEIHRD model, in Figures \ref{fig:Deterministic_vs_Stochastic_1}, \ref{fig:Deterministic_vs_Stochastic_2} and \ref{fig:Deterministic_vs_Stochastic_3} we present the solution of the deterministic model versus a number of possible \emph{scenarios} (simulations) obtained by the stochastic model. We can clearly observe that the outcomes by the stochastic model provide much more information about the future evolution of each model compartment, while the deterministic one seems to be somehow an \emph{averaged} version of the stochastic formulation. This gives just an initial insight of the potential of the newly introduced model, which will be completed in the next section with a more detailed analysis.
        
\section{Numerical and statistical analysis}
    
    
    The experiments have been conducted in a computer system with the following characteristics: CPU Intel Core i7-4720HQ 2.6GHz, 16GB RAM memory and GPU GeForce GTX 970M. The numerical codes have been implemented in Python programming language. We consider a equally spaced time grid, i.e. $\Delta t := t_{i+1} - t_i, \forall i$, with time step $\Delta t = \frac{1}{6}$ (around 4 hours).
    
    In this section, we perform a numerical and statistical study of the proposed stochastic model. As mentioned, the solution of the system of SDEs in \eqref{eq:sto_theta_seihrd} is a set of stochastic processes, meaning that we can ``extract'' a random variable at each time point. Doing so, not only a single value (like for the deterministic case), but also some statistics can be provided given a prescribed time. In this work, we consider the mean, the interquartile interval, $[Q_1, Q_3]$, with $Q_1$ and $Q_3$ being the first and the third quartiles, respectively, and a worst case scenario (WS), applied to both the evolution of the model variables and the possible model outputs.
    
    Here we take advantage of the experiments conducted in \cite{ivorra2020}, referred to the case of China. In Tables \ref{tab:params_inferred} and \ref{tab:params_calibrated} the reported values for the open coefficients are presented, distinguishing between the ones extracted from the literature or by experience, and the ones obtained by a calibration procedure in \cite{ivorra2020}, respectively. The initial data is given by $S(t_0) = N-1$, $E(t_0) = 1$ and $I(t_0) = I_u(t_0) = H_R(t_0) = H_D(t_0) = R_d(t_0) = R_u(t_0) = D(t_0) = 0$.
        \begin{table}[]
            \centering
            \begin{tabular}{c|c|c}
                \hline
                \textbf{Notation} & \textbf{Value} & \textbf{Description} \\
                \hline
                $N$ & $1400812636$ & Total population. \\
                $t_0$ & 1-12-2019 & Initial date.\\
                $T$ & 29-3-2020 & Final date. \\
                $\lambda_1$ & 23-1-2020 & Date when travel restrictions were imposed in Wuhan. \\
                $\lambda_2$ & 8-2-2020 & Inflexion date. \\
                $\underline{\theta}$ & $14\%$ & Percentage of documented cases at $\lambda_1$. \\
                $\overline{\theta}$ & $65\%$ & Percentage of documented cases at $\lambda_2$. \\
                $\alpha_H$ & $2.75\%$ & Percentage of infection produced by hospitalized people. \\
                $d_E$ & $5.5$ & Average days in compartment $E$. \\
                $d_I$ & $6.7$ & Average days in compartment $I$. \\
                $d_{I_u}$ & $14-d_I = 7.3$ & Average days in compartment $I_u$. \\
                $d_g$ & $6$ & Maximum reduction of $d_I$ due to the control measures. \\
                $C_o$ & $14$ & The period of convalescence. \\
                $p(t)$ & $1$ & Fraction of the infected people hospitalized. \\
                \hline
            \end{tabular}
            \caption{Parameters extracted from the experience and/or literature.}
            \label{tab:params_inferred}
        \end{table}
        
        \begin{table}[]
            \centering
            \begin{tabular}{c|c|c}
                \hline
                \textbf{Notation} & \textbf{Value} & \textbf{Description} \\
                \hline
                $\beta_I$ & $0.2887$ & Disease contact rate of a person in compartment $I$. \\
                $C_E$ & $0.3643$ & Reduction factor of the disease contact rate $\beta_E$ w.r.t $\beta_I$. \\
                $C_u$ & $0.4010$ & Reduction factor of the disease contact rate $\underline{\beta}_I$ w.r.t $\beta_I$. \\
                $\delta_R$ & $7.0000$ & Difference between days in compartment $H_R$ and $H_D$. \\
                $\delta_\omega$ & $0.0206$ & Difference between $\underline{\omega}$ and $\overline{\omega}$. \\
                $\underline{\omega}$ & $0.0157$ & Lower bound of the fatality rate. \\
                $\kappa_1$ & $0.1082$ & Efficiency of the control measures. \\
                \hline
            \end{tabular}
            \caption{Parameters obtained by calibration to the data.}
            \label{tab:params_calibrated}
        \end{table}

    \subsection{Model variables}
    
        We firstly test the evolution of the model variables. Thus, we extract some simulation-based statistics at a couple of time instants, the 8th February (inflection point) and the 29th March (final point). Furthermore, the impact of different levels of uncertainty is also reported by considering several representative values for $\sigma_{\beta_I}$, reflecting situations of no uncertainty ($\sigma_{\beta_I} = 0$), low uncertainty ($\sigma_{\beta_I} = 0.1$) and high uncertainty ($\sigma_{\beta_I} = 0.5$). In all cases, the long-term mean of the perturbation is set to the calibrated value of $\beta_I$ for the deterministic model, i.e. $\mu_{\beta_I} = \beta_I$. The mean reverting speed parameter is chosen as $\nu_{\beta_I} = 1$, thus representing a regular (not too high, not too low) reversion speed. In Table \ref{tab:variables_CIR}, the obtained results are presented. We can clearly observe the significant impact of the uncertainty in the disease evolution. In the case of higher uncertainty, the number of infections, in average, is almost doubled and, the worst case scenario multiplies this value by six. Even when a lower volatility is considered, the increment of cases and deaths in the worst scenario becomes important, up to $50\%$.

            \begin{table}[h]
                \centering
                \begin{tabular}{c||c|c c c|c c c}
                    \hline
                    & \multicolumn{7}{c}{8th February, 2020 ($t=69$)} \\
                    \hline
                    & $\sigma_{\beta_I} = 0$ & \multicolumn{3}{c|}{$\sigma_{\beta_I} = 0.1$} & \multicolumn{3}{c}{$\sigma_{\beta_I} = 0.5$} \\
                    \hline
                    & Mean & Mean & $[Q_1, Q_3]$ & WS ($95\%$) & Mean & $[Q_1, Q_3]$ & WS ($95\%$) \\
                    \hline\hline
                    $E(t)$ & $2993$ & $3067$ & $[2506, 3519]$ & $4510$ & $5401$ & $[1049, 5415]$ & $18970$ \\
                    $I(t)$ & $1340$ & $1376$ & $[1125, 1578]$ & $2017$ & $2419$ & $[476, 2423]$ & $8522$ \\
                    $I_u(t)$ & $3854$ & $3945$ & $[3249, 4505]$ & $5724$ & $6811$ & $[1434, 6940]$ & $23728$ \\
                    $H_R(t)$ & $3252$ & $3328$ & $[2732, 3806]$ & $4854$ & $5799$ & $[1182, 5863]$ & $20340$ \\
                    $H_D(t)$ & $214$ & $219$ & $[181, 250]$ & $318$ & $377$ & $[80, 386]$ & $1311$ \\
                    $R_d(t)$ & $1846$ & $1888$ & $[1559, 2153]$ & $2726$ & $3231$ & $[701, 3317]$ & $11168$ \\
                    $R_u(t)$ & $4296$ & $4390$ & $[3656, 4985]$ & $6238$ & $7301$ & $[1738, 7690]$ & $24654$ \\
                    $D(t)$ & $131$ & $134$ & $[112, 152]$ & $190$ & $222$ & $[53, 235]$ & $747$ \\
                    \hline
                    \hline
                    \hline
                    & \multicolumn{7}{c}{29th March, 2020 ($t=119$)} \\
                    \hline
                    & $\sigma_{\beta_I} = 0$ & \multicolumn{3}{c|}{$\sigma_{\beta_I} = 0.1$} & \multicolumn{3}{c}{$\sigma_{\beta_I} = 0.5$} \\
                    \hline
                    & Mean & Mean & $[Q_1, Q_3]$ & WS ($95\%$) & Mean & $[Q_1, Q_3]$ & WS ($95\%$) \\
                    \hline\hline
                    $E(t)$ & $1$ & $1$ & $[1, 1]$ & $2$ & $2$ & $[0, 3]$ & $10$ \\
                    $I(t)$ & $0$ & $0$ & $[0, 0]$ & $0$ & $0$ & $[0, 0]$ & $1$ \\
                    $I_u(t)$ & $173$ & $177$ & $[146, 203]$ & $259$ & $309$ & $[63, 314]$ & $1075$ \\
                    $H_R(t)$ & $232$ & $237$ & $[194, 272]$ & $348$ & $416$ & $[83, 420]$ & $1458$ \\
                    $H_D(t)$ & $30$ & $30$ & $[25, 35]$ & $45$ & $53$ & $[11, 54]$ & $186$ \\
                    $R_d(t)$ & $8460$ & $8662$ & $[7118, 9910]$ & $12624$ & $15087$ & $[3101, 15287]$ & $52651$ \\
                    $R_u(t)$ & $9969$ & $10198$ & $[8442, 11616]$ & $14681$ & $17386$ & $[3862, 17941]$ & $59616$ \\
                    $D(t)$ & $417$ & $426$ & $[353, 486]$ & $614$ & $728$ & $[161, 751]$ & $2502$ \\
                    \hline
                \end{tabular}
                \caption{Variables of the Stochastic $\theta$-SEIHRD model: $\nu_{\beta_I} = 1$, $\mu_{\beta_I} = \beta_I$ and $n = 2^{15}$ Monte Carlo simulations. Columns: mean, interquartile interval ($[Q_1, Q_3]$) and worst case scenario (WS).}
                \label{tab:variables_CIR}
            \end{table}
            
\begin{figure}[h!]
                \centering
                    \subfigure[$\sigma_{\beta_I} = 0.1$]{\includegraphics[width=0.49\textwidth]{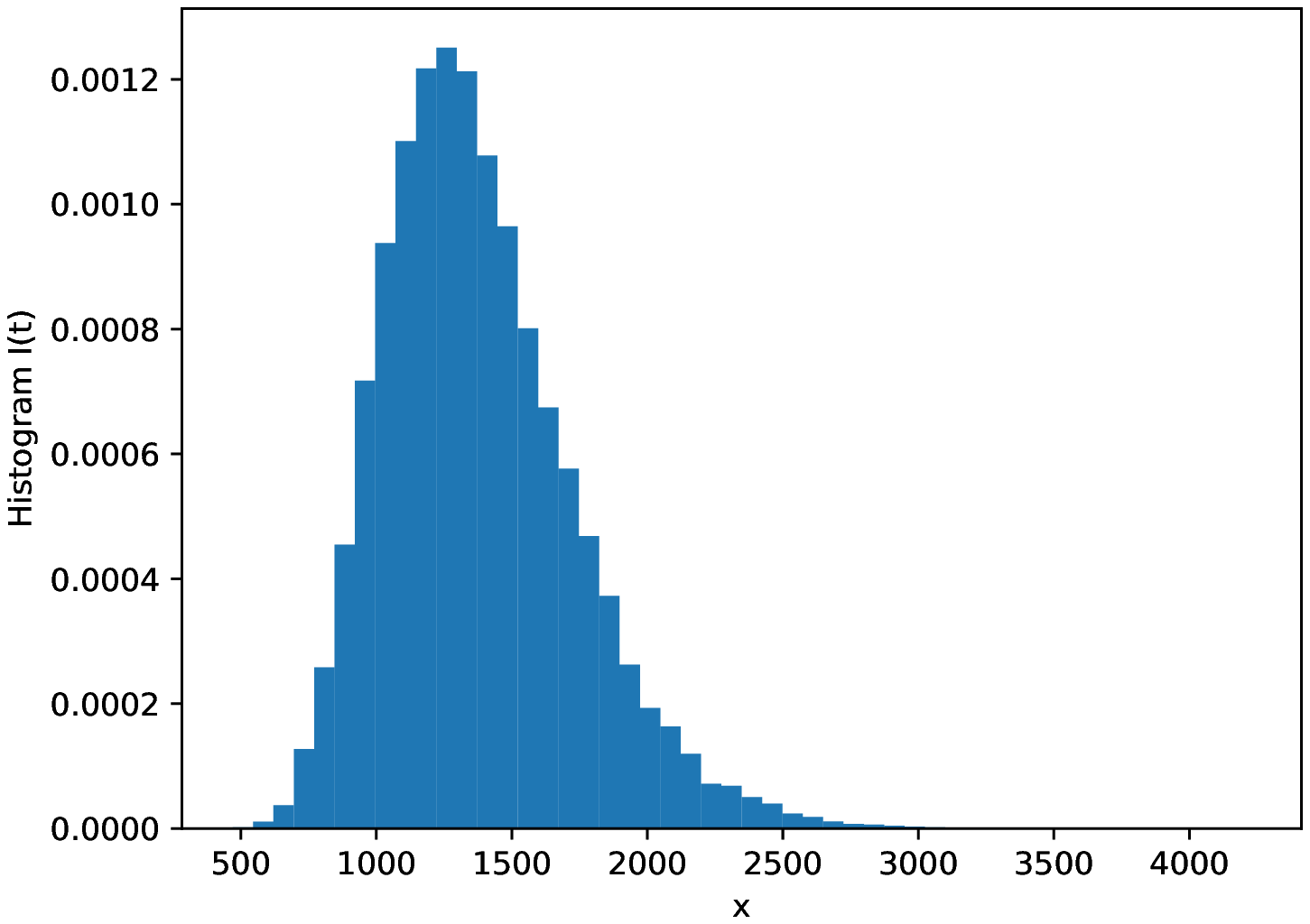}}
                    \subfigure[$\sigma_{\beta_I} = 0.5$]{\includegraphics[width=0.49\textwidth]{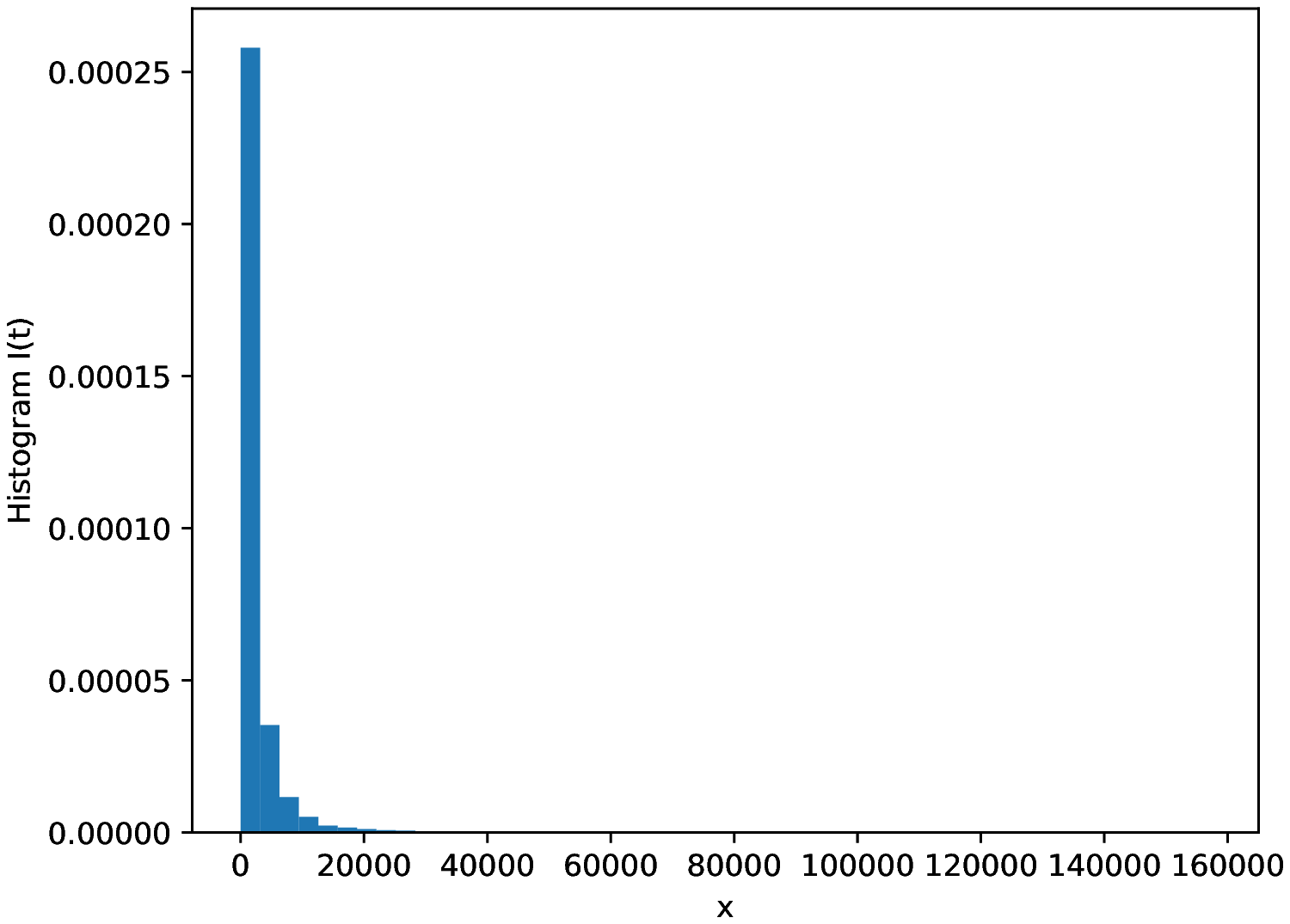}}
                \caption{Histogram of $I(t)$. Setting: $\nu_{\beta_I} = 1$ and $\mu_{\beta_I} = \beta_I$, with $n = 2^{15}$ Monte Carlo simulations.}
                \label{fig:histograms_I}
            \end{figure}

        Next to the previous experiment, in Figure \ref{fig:histograms_I} we present the histograms for the model variable $I(t)$ to give an insight of the impact of the uncertainty in the infection evolution. First, we clearly observe that the produced distribution is skewed with a fatter right tail. Secondly, the bigger the volatility of process $\tilde{\beta}_I(t)$, i.e. $\sigma_{\beta_I}$, is, the fatter the right tail becomes, thus indicating more probability of extreme events. 

    \subsection{Outputs}
    
        In \cite{ivorra2020}, the authors proposed a set of possible outputs that can be useful for the authorities to plan the resources allocation, like the number of clinical beds, for example, among other indicators. These outputs are:
        \begin{itemize}
            \item The cumulative number of COVID-19 cases, $c_m(t)$, at time $t$:
                \begin{equation*}
                    c_m(t) = H_R(t) + H_D(t) + R_d(t) + D(t) = c_m(t_0) + \int_{t_0}^t \theta(s)\gamma_I(s)I(s)\ds.
                \end{equation*}
            
            \item The cumulative number of deaths due to the COVID-19, at time $t$: $d_m(t) = D(t)$.
            
            \item The basic reproduction number, $R_0$, and the effective reproduction number, $R_e(t)$, at time $t$, where $R_0 = R_e(t_0)$, and
                \begin{equation*}
                    R_e(t) = \frac{U_e(t)}{\gamma_E\gamma_I(t)\gamma_{H_R}(t)\gamma_{H_D}(t)\gamma_{I_u}(t)}\frac{S(t)}{N},
                \end{equation*}
                with\footnote{The dependence of the coefficients on the time $t$ has been omitted for notational purposes. All the coefficients take their particular values at day $t$.}
                \begin{equation*}
                \begin{aligned}
                    U_e(t) &= \left(\left((m_{I_u}\beta_{I_u}(1 - \theta)\gamma_{H_R} + m_{H_R}\beta_{H_R}\gamma_{I_u}(\theta - \omega))\gamma_I + m_I\beta_I\gamma_{H_R}\gamma_{I_u}\right)\gamma_E + m_E\beta_E\gamma_I\gamma_{H_R}\gamma_{I_u} \right)\gamma_{H_D} \\
                    &+ m_{H_D}\beta_{H_D}\omega\gamma_E\gamma_I\gamma_{H_R}\gamma_{I_u}.
                \end{aligned}
                \end{equation*}
                
            \item Hospitalized people, $\textrm{Hos}(t)$, at time $t$:
                \begin{equation*}
                    \textrm{Hos}(t) = H_D(t) + p(t)\left( H_R(t) + R_d(t) - R_d(t - C_o) \right),
                \end{equation*}
                where $p(t)$ is the fraction of people in compartment $H_R$ that are hospitalized and $C_o$ is the period of convalescence.
                
            \item Maximum number of hospitalized people in the interval $[t_0, t]$:
                \begin{equation*}
                    \textrm{MHos}(t) = \max_{\tau\in[t_0, t]} \textrm{Hos}(\tau).
                \end{equation*}
            
            \item The number of individuals infected by others belonging to compartments $E$, $I_u$ and $H = H_R + H_D$:
                \begin{equation*}
                \begin{aligned}
                    \Gamma_E(t) &= \int_{t_0}^t m_E(s)\beta_E E(s)\frac{S(s)}{N}\ds, \\
                    \Gamma_{I_u}(t) &= \int_{t_0}^t m_{I_u}(s)\beta_{I_u} I_u(s)\frac{S(s)}{N}\ds, \\
                    \Gamma_H(t) &= \int_{t_0}^t (m_{H_R}(s)\beta_{H_R} H_R(s) + m_{H_D}(s)\beta_{H_D} H_D(s))\frac{S(s)}{N}\ds,
                \end{aligned}
                \end{equation*}
                respectively.
        \end{itemize}
            
        We therefore perform a similar statistical analysis as before, although now reporting the model outputs. The results are shown in Table \ref{tab:outputs_CIR}. Again, it is clear that the uncertainty can significantly affect the disease evolution, and it should be taken into consideration. For example, the people requiring hospitalization may vary from around $4500$ with no stochasticity included, to around $30000$ including randomness.

            \begin{table}[h]
                \centering
                \begin{tabular}{c||c|c c c|c c c}
                    \hline
                    & \multicolumn{7}{c}{8th February, 2020 ($t=69$)} \\
                    \hline
                    & $\sigma_{\beta_I} = 0$ & \multicolumn{3}{c|}{$\sigma_{\beta_I} = 0.1$} & \multicolumn{3}{c}{$\sigma_{\beta_I} = 0.5$} \\
                    \hline
                    & Mean & Mean & $[Q_1, Q_3]$ & WS ($95\%$) & Mean & $[Q_1, Q_3]$ & WS ($95\%$) \\
                    \hline\hline
                    $c_m(t)$ & $5440$ & $5571$ & $[4586, 6362]$ & $8088$ & $9631$ & $[2026, 9812]$ & $33466$ \\
                    $d_m(t)$ & $131$ & $134$ & $[112, 152]$ & $190$ & $222$ & $[53, 235]$ & $747$ \\
                    $R_e(t)$ & $0.3363$ & $0.3364$ & $[0.3151$, $0.3562]$ & $0.3891$ & $0.3367$ & $[0.2245$, $0.4149]$ & $0.6340$ \\
                    $Hos(t)$ & $4040$ & $4134$ & $[3395, 4727]$ & $6026$ & $7197$ & $[1471, 7273]$ & $25262$ \\
                    $\text{MHos}(t)$ & $4040$ & $4134$ & $[3395, 4727]$ & $6026$ & $7197$ & $[1471, 7273]$ & $25262$ \\
                    $\Gamma_E(t)$ & $5012$ & $5126$ & $[4255, 5833]$ & $7337$ & $8625$ & $[1979, 9013]$ & $29414$ \\
                    $\Gamma_{I_u}(t)$ & $4550$ & $4646$ & $[3864, 5285]$ & $6640$ & $7600$ & $[1764, 7976]$ & $25755$ \\
                    $\Gamma_H(t)$ & $198$ & $202$ & $[168, 230]$ & $288$ & $328$ & $[78, 346]$ & $1106$ \\
                    \hline 
                    \hline
                    \hline
                    & \multicolumn{7}{c}{29th March, 2020 ($t=119$)} \\
                    \hline
                    & $\sigma_{\beta_I} = 0$ & \multicolumn{3}{c|}{$\sigma_{\beta_I} = 0.1$} & \multicolumn{3}{c}{$\sigma_{\beta_I} = 0.5$} \\
                    \hline
                    & Mean & Mean & $[Q_1, Q_3]$ & WS ($95\%$) & Mean & $[Q_1, Q_3]$ & WS ($95\%$) \\
                    \hline\hline
                    $c_m(t)$ & $9140$ & $9358$ & $[7691, 10704]$ & $13631$ & $16286$ & $[3358, 16526]$ & $56752$ \\
                    $d_m(t)$ & $417$ & $426$ & $[353, 486]$ & $614$ & $728$ & $[161, 751]$ & $2502$ \\
                    $R_e(t)$ & $0.0013$ & $0.0013$ & $[0.0012$ , $0.0014]$& $0.0015$ & $0.0013$ & $[0.0009$, $0.0016]$ & $0.0024$ \\
                    $\text{Hos}(t)$ & $306$ & $314$ & $[257, 360]$ & $459$ & $549$ & $[111, 555]$ & $1927$ \\
                    $\text{MHos}(t)$ & $4558$ & $4671$ & $[3832, 5347]$ & $6816$ & $8195$ & $[1662, 8258]$ & $28681$ \\
                    $\Gamma_E(t)$ & $5259$ & $5379$ & $[4464, 6122]$ & $7705$ & $9070$ & $[2073, 9465]$ & $30925$ \\
                    $\Gamma_{I_u}(t)$ & $5388$ & $5504$ & $[4570, 6264]$ & $7886$ & $9082$ & $[2080, 9479]$ & $30911$ \\
                    $\Gamma_H(t)$ & $229$ & $234$ & $[195, 266]$ & $334$ & $384$ & $[89, 402]$ & $1298$ \\
                    \hline
                \end{tabular}
                \caption{Outputs of the Stochastic $\theta$-SEIHRD model: $\nu_{\beta_I} = 1$, $\mu_{\beta_I} = \beta_I$ and $n = 2^{15}$ Monte Carlo simulations. Columns: mean, interquartile interval ($[Q_1, Q_3]$) and worst case scenario (WS).}
                \label{tab:outputs_CIR}
            \end{table}

\begin{figure}[h!]
                \centering
                    \subfigure[Reported cumulative]{\includegraphics[width=0.49\textwidth]{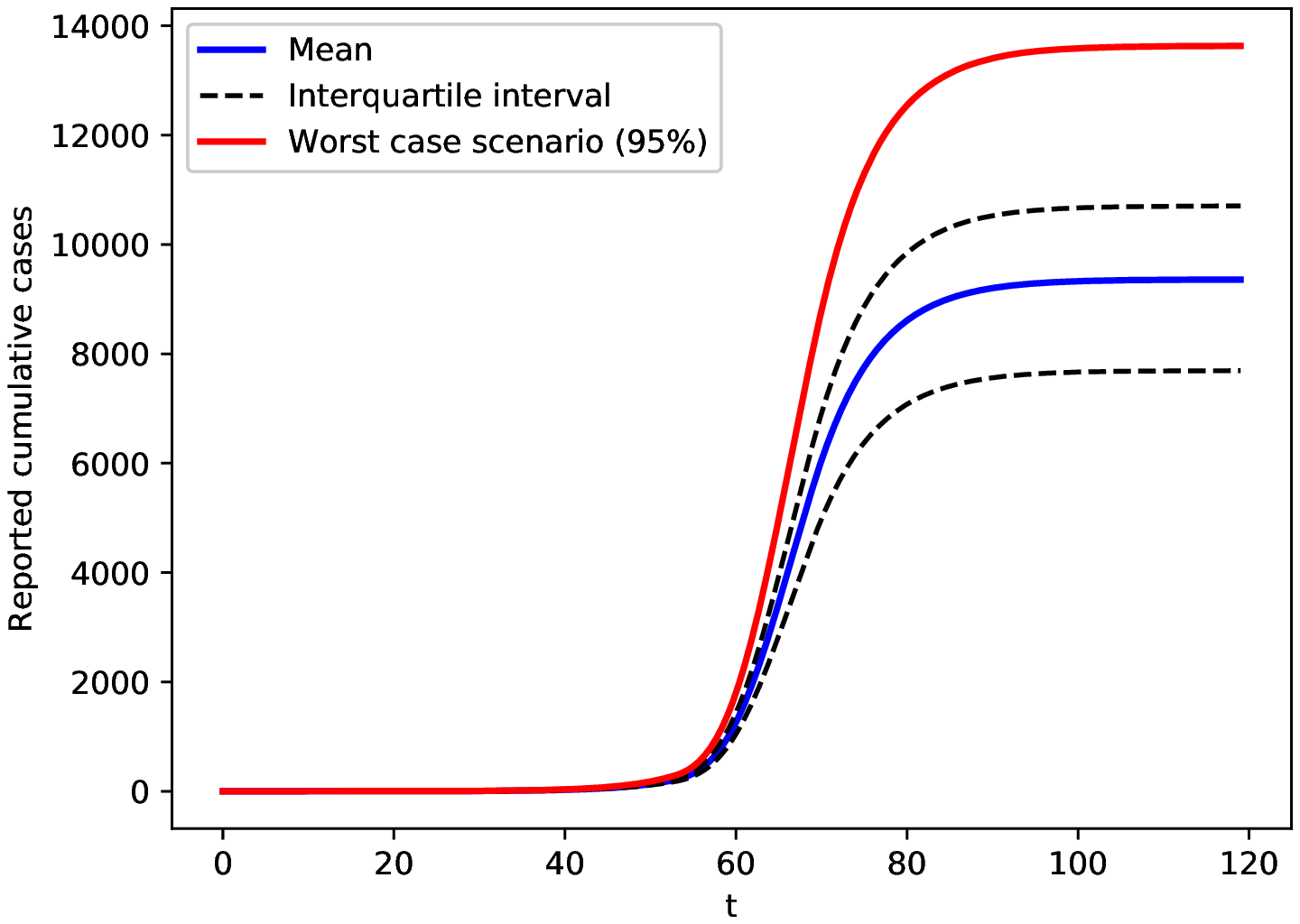}}
                    \subfigure[Undetected cumulative]{\includegraphics[width=0.49\textwidth]{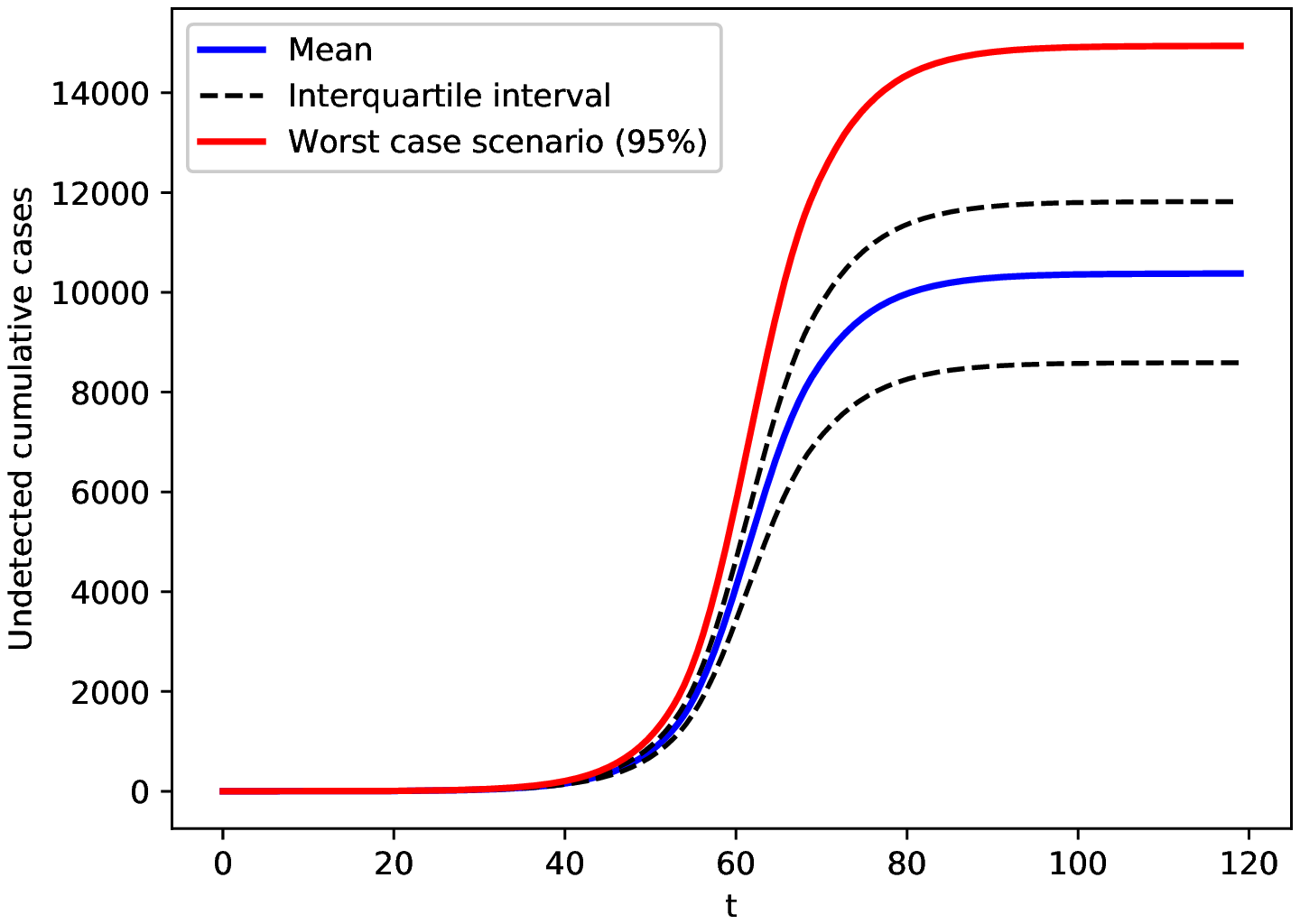}}
                    \subfigure[Hospitalized]{\includegraphics[width=0.49\textwidth]{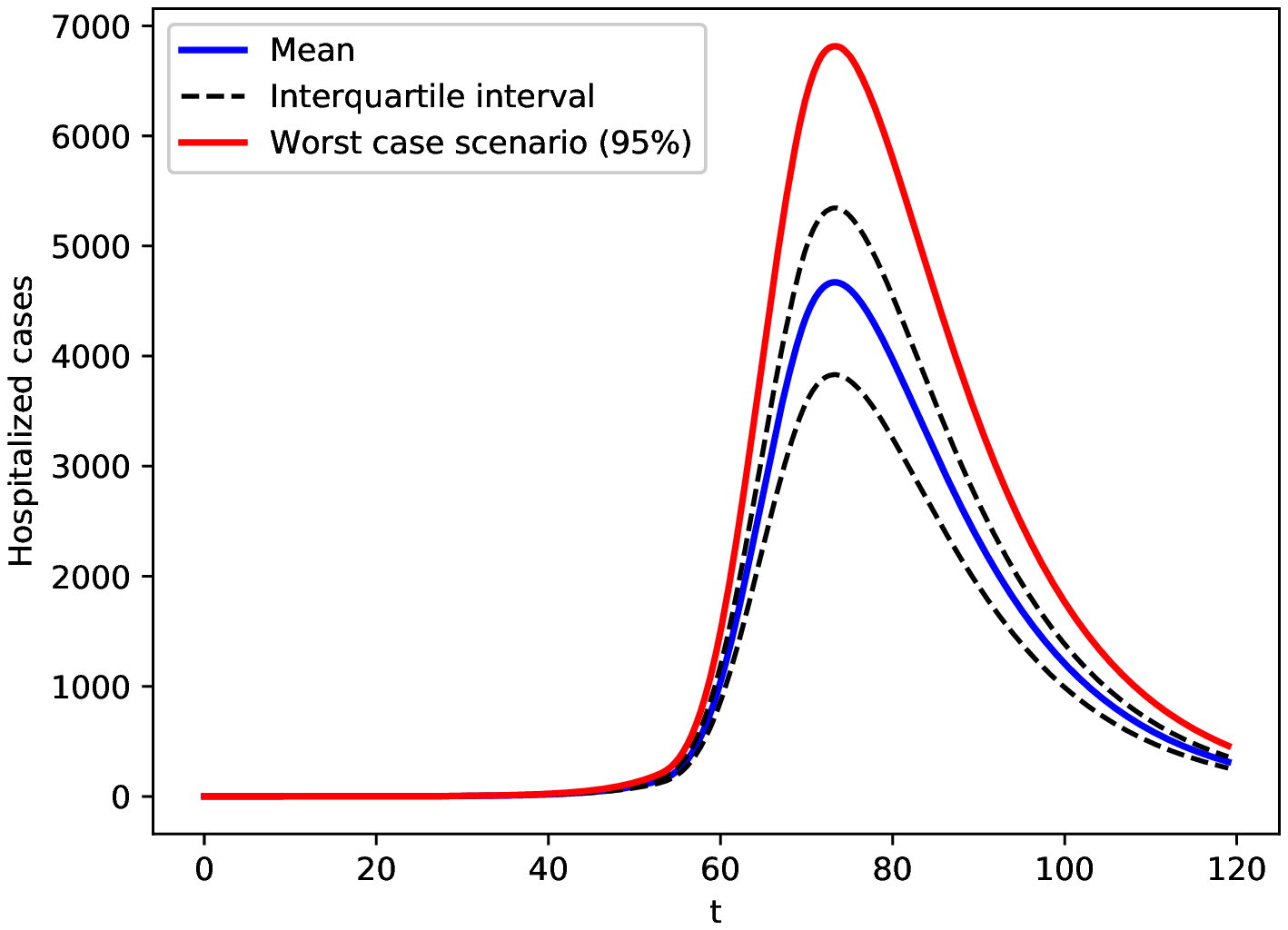}}
                    \subfigure[Reported daily infected]{\includegraphics[width=0.49\textwidth]{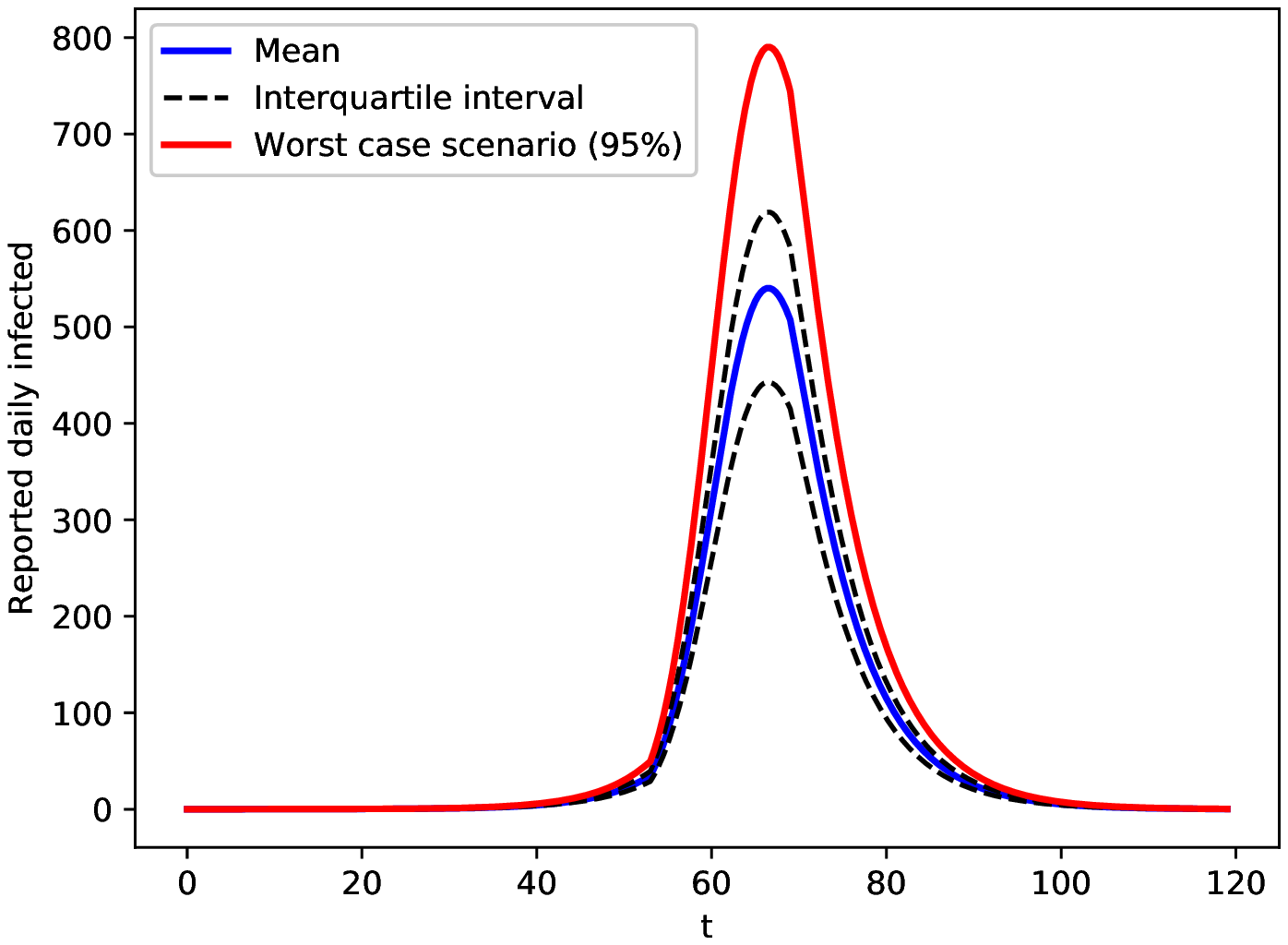}}
                    \subfigure[Reported daily deaths]{\includegraphics[width=0.49\textwidth]{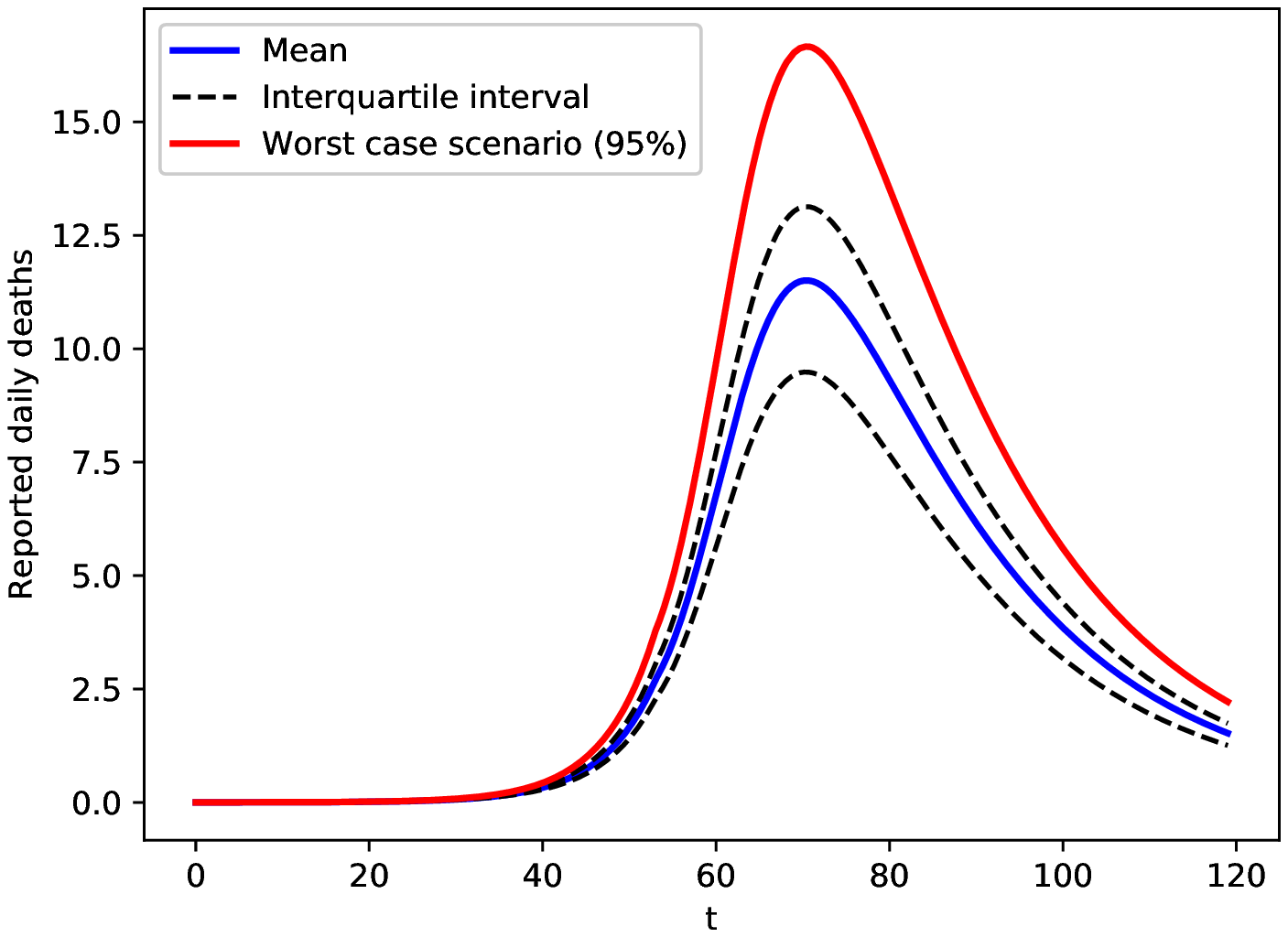}}
                    \subfigure[Reported daily recovered]{\includegraphics[width=0.49\textwidth]{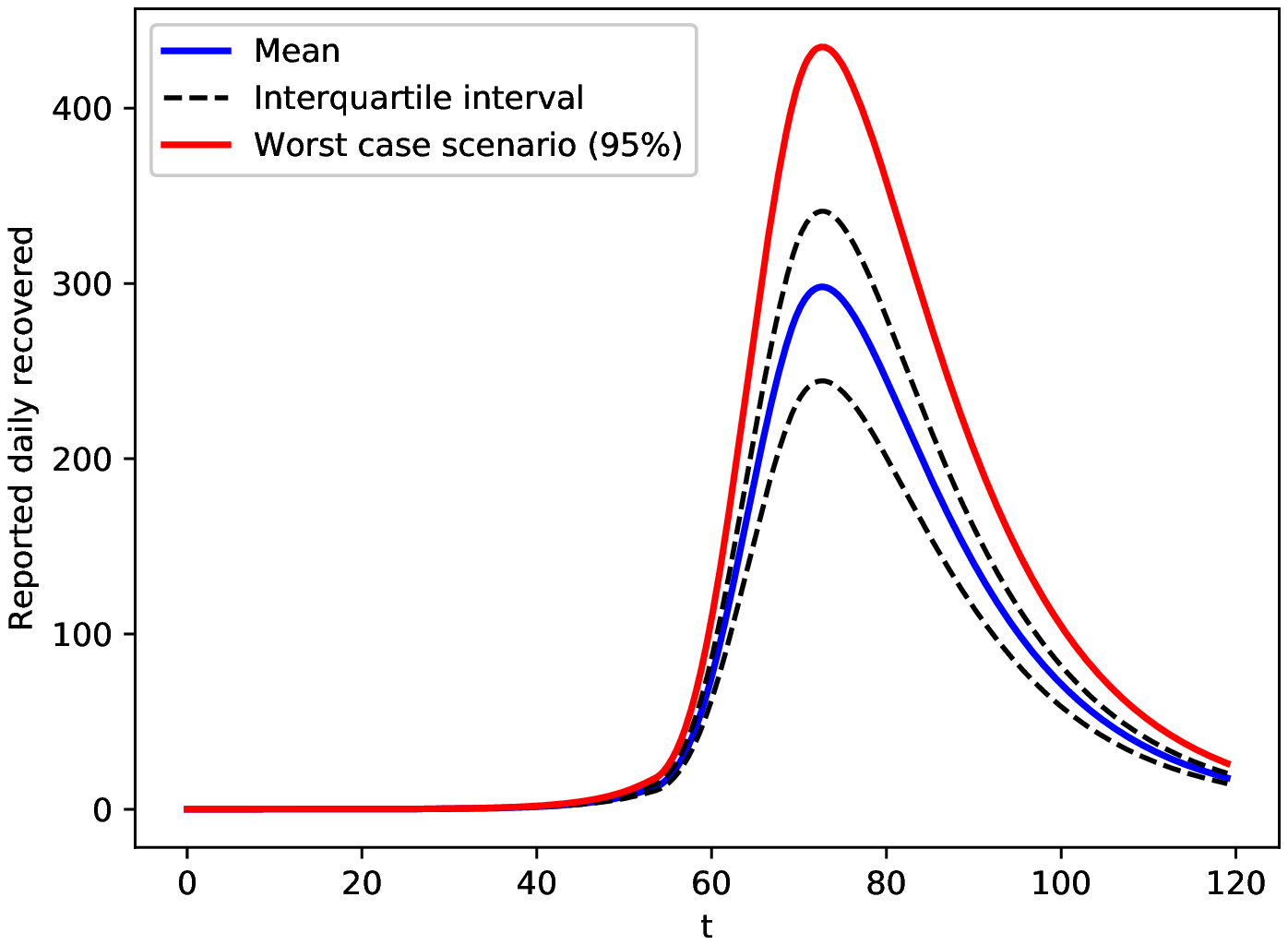}}
                \caption{Epidemic curves: mean, interquartile interval and worst case scenario. Setting: $\nu_{\beta_I} = 1$, $\mu_{\beta_I} = \beta_I$ and $\sigma_{\beta_I} = 0.1$, with $n = 2^{15}$ Monte Carlo simulations.}
                \label{fig:epidemic_curves_01}
            \end{figure}

\begin{figure}[h!]
                \centering
                    \subfigure[Reported cumulative]{\includegraphics[width=0.49\textwidth]{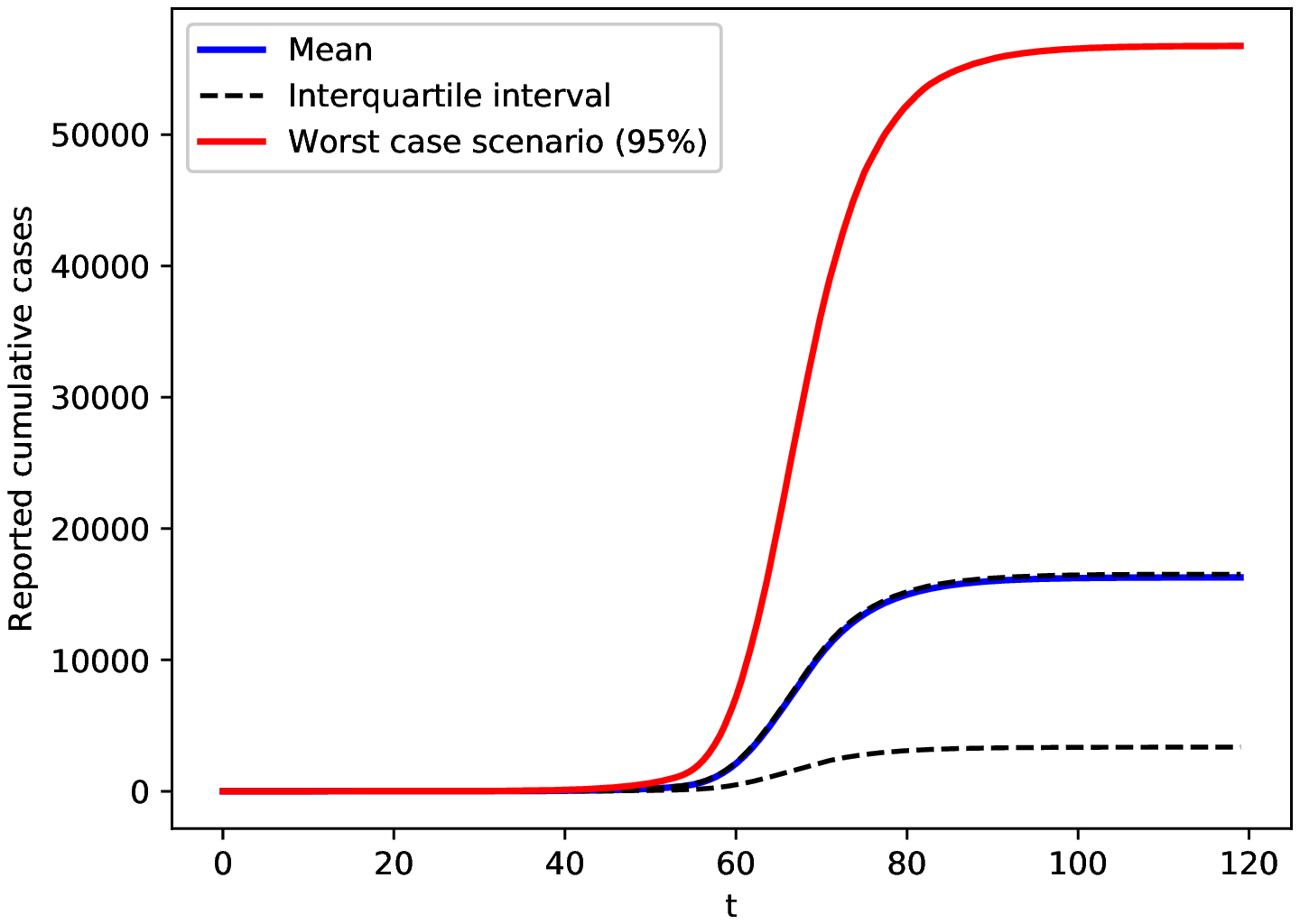}}
                    \subfigure[Undetected cumulative]{\includegraphics[width=0.49\textwidth]{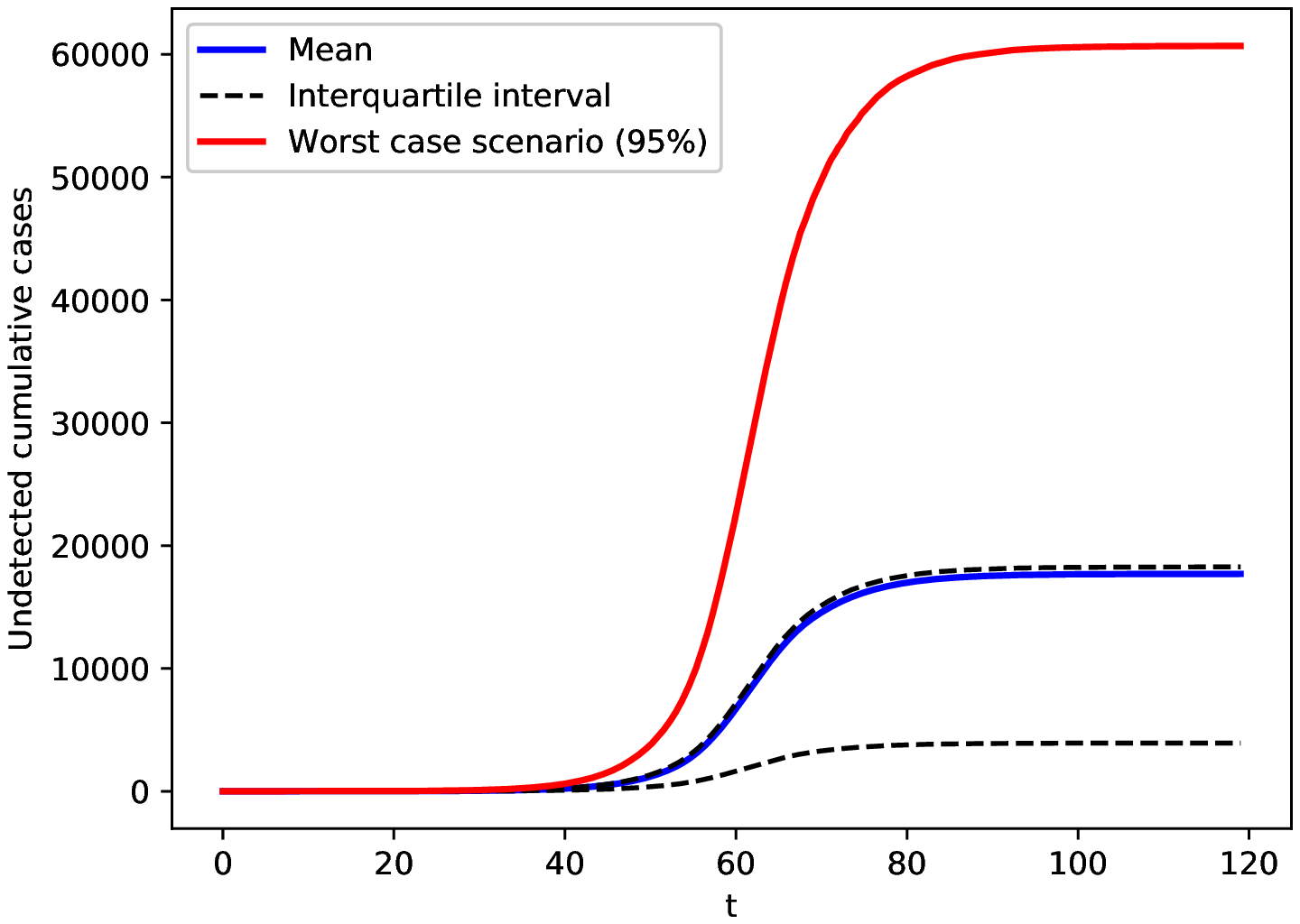}}
                    \subfigure[Hospitalized]{\includegraphics[width=0.49\textwidth]{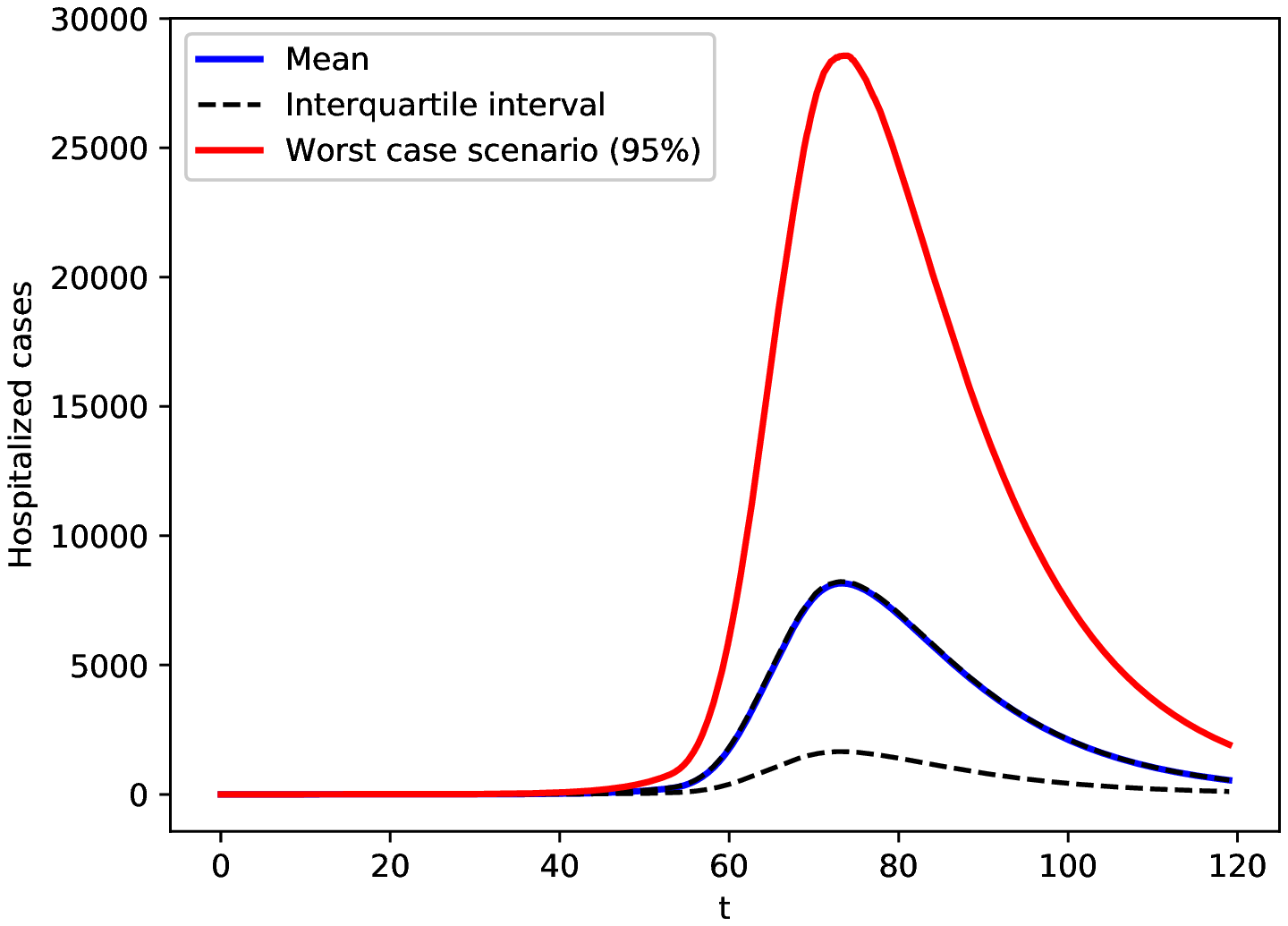}}
                    \subfigure[Reported daily infected]{\includegraphics[width=0.49\textwidth]{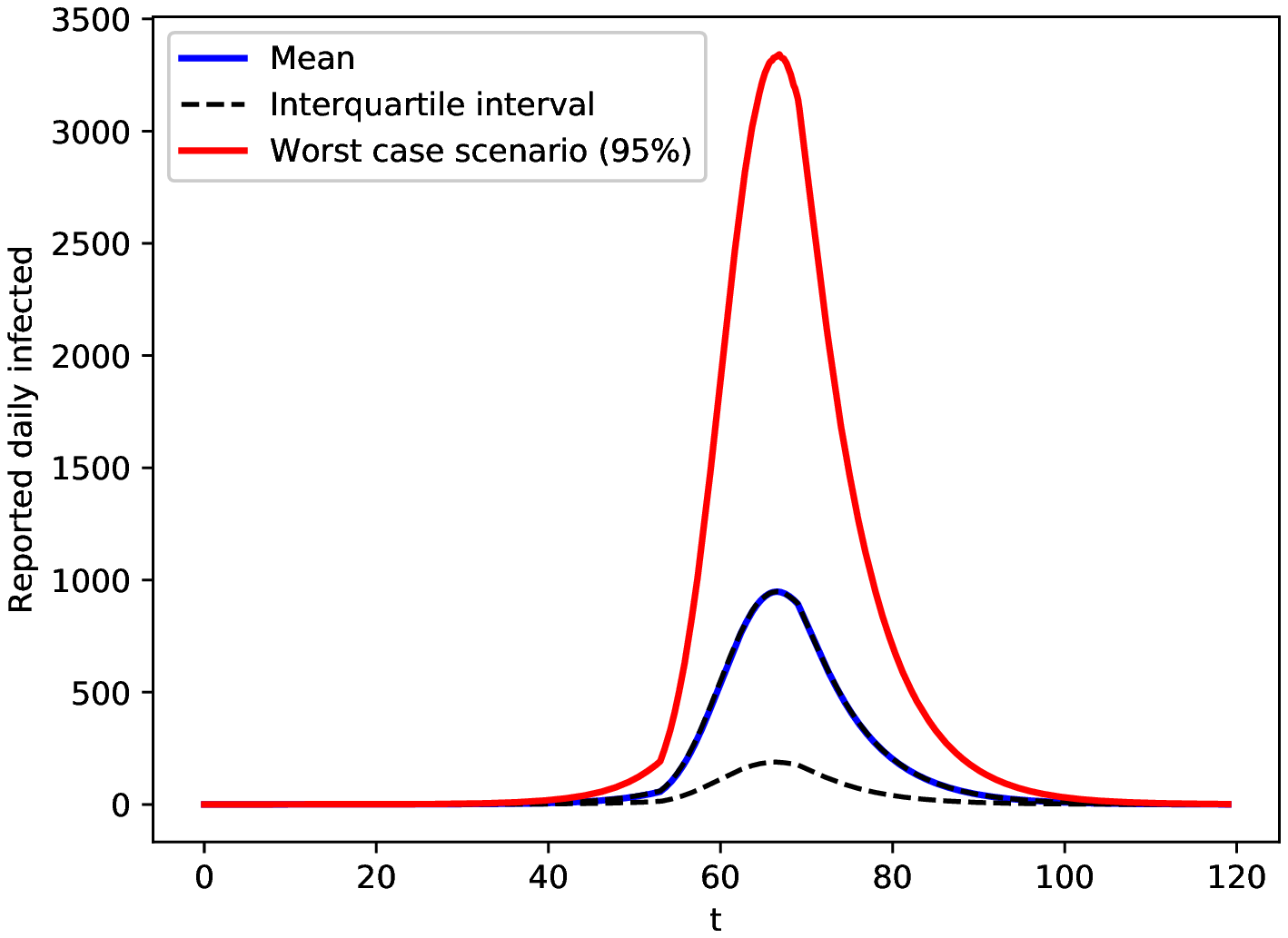}}
                    \subfigure[Reported daily deaths]{\includegraphics[width=0.49\textwidth]{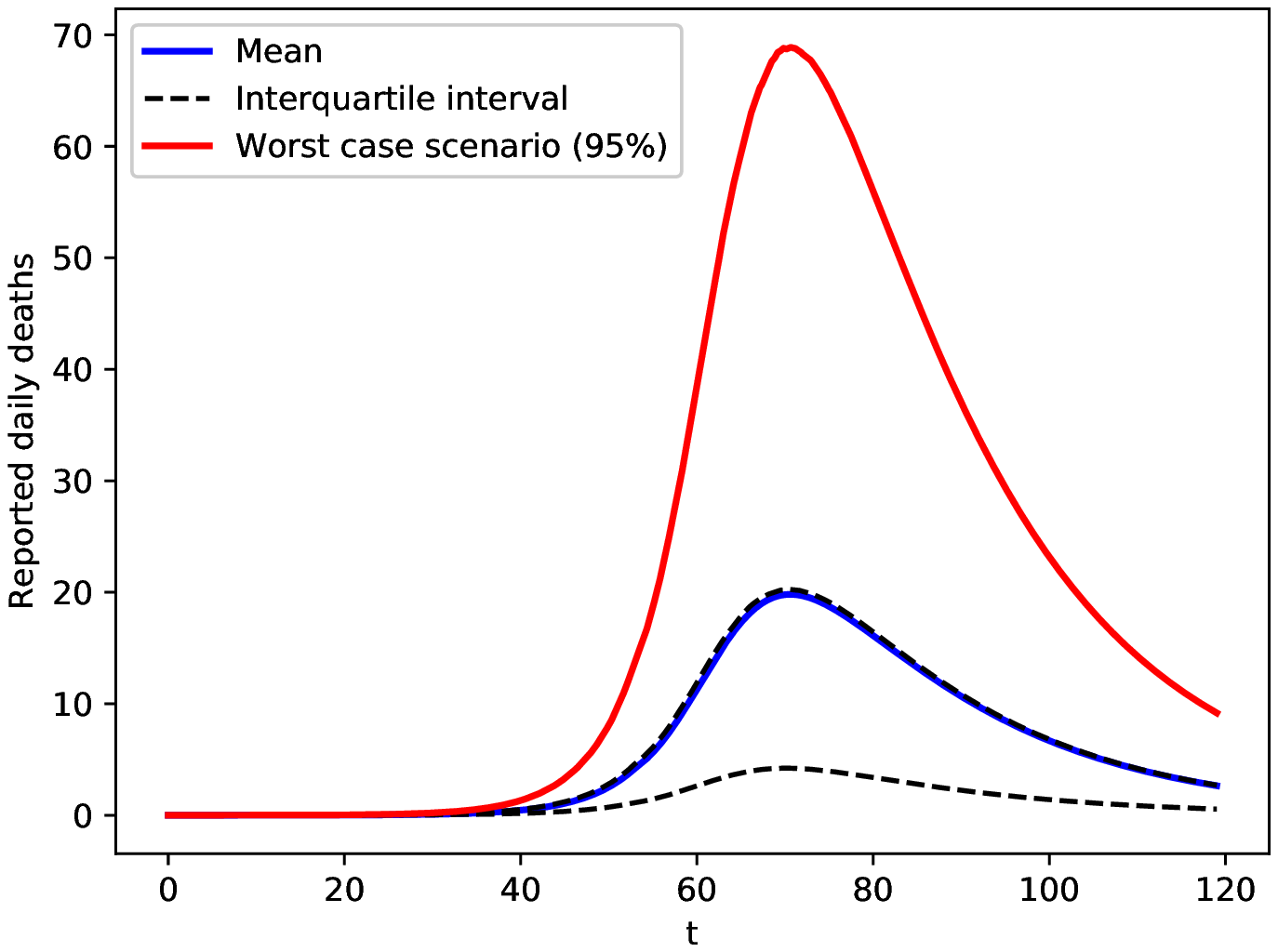}}
                    \subfigure[Reported daily recovered]{\includegraphics[width=0.49\textwidth]{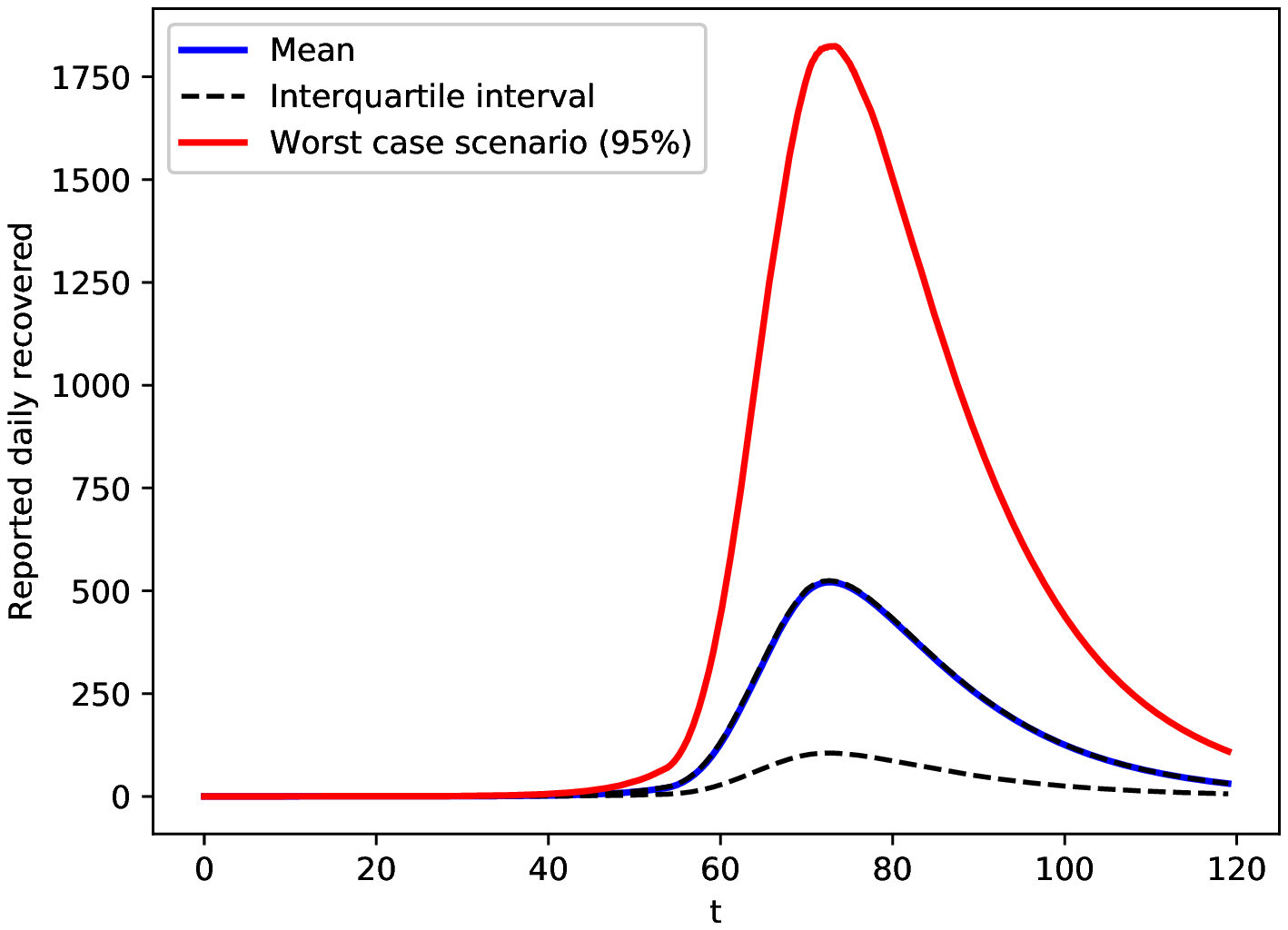}}
                \caption{Epidemic curves: mean, interquartile interval and worst case scenario. Setting: $\nu_{\beta_I} = 1$, $\mu_{\beta_I} = \beta_I$ and $\sigma_{\beta_I} = 0.5$, with $n = 2^{15}$ Monte Carlo simulations.}
                \label{fig:epidemic_curves_05}
            \end{figure}

        As previously pointed out, one of the major advantages of the $\theta$-SEIHRD model is that it accounts for important aspects of the COVID-19 pandemic, directly affecting the population or the healthcare systems. Particularly, the evolution of some curves, infected, hospitalized and deaths, is typically reported by the authorities, in both cumulative and daily fashion. In Figures \ref{fig:epidemic_curves_01} and \ref{fig:epidemic_curves_05}, we show the model outcomes for these specific curves, considering two uncertainty levels, $\sigma_{\beta_I} = 0.1$ and $\sigma_{\beta_I} = 0.5$, respectively. Again, we present the mean, the interquartile interval and the worst case scenario. From this experiment, an interesting observation can be extracted. Looking at the different patterns in the restults w.r.t. the volatility parameter, we can see that an increasing uncertainty pushes the mean close to the third quartile, $Q_3$, meaning that the disease evolves, in average, according to the $75\%$ worst scenario. This fact gives an insight of how important the randomness can be and why it is crucial to include it in the modelling.
            
            

\section{Discussion and conclusions}\label{sec:discussion}

    We have extended the model developed in \cite{ivorra2020} by incorporating randomness to some relevant coefficients and we have shown the importance of considering this uncertainty. In this way, besides the information provided by the deterministic model (which allows to obtain a proxy to the average of the main variables), we can take advantage of a more complete modelling approach, which allows not only to compute confidence intervals for these variables in this new random setting but also to obtain the worst case scenario. The information about the model variables in this worst case scenario allows to develop more conservative policies in the actions against the results of the COVID-19, as for example to plan larger health resources to take care of a larger number of people requiring hospitalization at different levels. The differences between the deterministic and stochastic models in terms of the information contained in the output variables has been clearly illustrated in the previous section.

    However, we also understand that research in this line can be extended. A fist possible extension comes from making the parameters independent among each other and use different stochastic processes to characterize their dynamics. In this first approach, as in \cite{ivorra2020}, we consider that all parameters depend on $\beta_I$. Also, it seems possible to incorporate randomness to the more recent model in \cite{ivorra2020b} that mainly incorporate a new compartment of persons in quarantine (Q) to model the situation in certain countries like Italy.

\section*{Acknowledgements}
    
    This work was funded by Xunta de Galicia grant ED431C2018/033, including FEDER funding. Both authors also acknowledge the support received from the Centro de Investigaci\'on de Galicia ``CITIC'', funded by Xunta de Galicia and the European Union (European Regional Development Fund- Galicia 2014-2020 Program), by grant ED431G 2019/01.
    
    \'Alvaro Leitao acknowledges the financial support from the Spanish Ministry of Science, Innovation and Universities, through the Juan de la Cierva-formaci\'on 2017 (FJC17) grant in the framework of the national programme for R\&D 2013-2016.

\bibliographystyle{plain}
\bibliography{stochastic_theta_SEIHRD}

\end{document}